\documentclass[11pt]{article}
\usepackage[left=1in,right=1in,top=1in,bottom=1in]{geometry}
\usepackage{amsmath}
\usepackage{amssymb}
\usepackage{mathtools}
\usepackage{bm}
\usepackage{graphicx,float}
\usepackage{tikz}
\usepackage[title]{appendix}

\title{A gyroscopic polynomial basis in the sphere}
\author{Abram C. Ellison, Keith Julien, Geoffrey M. Vasil}

\makeatletter
\newcommand*\bigcdot{\mathpalette\bigcdot@{.6}}
\newcommand*\bigcdot@[2]{\mathbin{\vcenter{\hbox{\scalebox{#2}{$\m@th#1\bullet$}}}}}
\makeatother
\newcommand\floor[1]{\left\lfloor #1 \right\rfloor}

\renewcommand{\vec}{\boldsymbol}
\newcommand{\del}{\vec{\nabla}}
\newcommand{\norm}[1]{\left\lVert#1\right\rVert}
\newcommand{\dd}{\,\text{d}}

\usepackage{xcolor}

\newcommand{\zero}{0\hspace{0.5ex}}

\begin{document}

\maketitle

\begin{abstract}
Standard spectral codes for full sphere dynamics utilize a combination of spherical harmonics and a suitable
radial basis to represent fluid variables.  These basis functions have a rotational invariance not present in
geophysical flows.  Gyroscopic alignment - alignment of dynamics along the axis of rotation - is a
hallmark of geophysical fluids in the rapidly rotating regime.  The Taylor-Proudman theorem, resulting
from a dominant balance of the Coriolis force and the pressure gradient force, yields nearly invariant flows along this axial direction.
In this paper we tailor a coordinate system to the cylindrical structures found in rotating spherical flows.
This ``spherindrical'' coordinate system yields a natural hierarchy of basis functions, composed of Jacobi polynomials
in the radial and vertical direction, regular throughout the ball.
We expand fluid variables using this basis and utilize sparse Jacobi polynomial algebra to implement all operators
relevant for partial differential equations in the spherical setting.  We demonstrate the representation power of
the basis in three eigenvalue problems for rotating fluids.

\end{abstract}

\section{Introduction}
Spherical geometry is a critical setting for three-dimensional physics simulations.  The geo- and astrophysics community
require efficient approaches to simulating fluid flow in the ball environment.  Native parameter regimes for these applications
are extreme, leading to computationally demanding simulations requiring prohibitive spatio-temporal resolutions.  Rotational constraint is common to many of these fluid
dynamics investigations.  This paper presents a sparse numerical approach in spherical geometry that leverages the gyroscopic alignment
so prevalent in celestial bodies.

Spectral codes for the ball utilize an expansion in modes to represent fields \cite{Fornberg_1998, Trefethen_2000, Boyd_2001, Burns_Vasil_Oishi_Lecoanet_Brown_2020}.
For spherical geometries, the coordinate system in which we work strongly influences our choice of basis for expansion.   The natural choice is to utilize spherical polar
coordinates ${ \left( r, \theta, \phi \right) }$ to denote position in the ball, where ${ 0 \le r \le 1 }$, ${ 0 \le \theta \le \pi }$ and
${ 0 \le \phi < 2 \pi }$.  This choice immediately implies the
Fourier basis in azimuthal angle $\phi$ due to its periodicity.  To be concrete, we represent a field $f(r, \theta, \phi)$
as
\begin{equation}
f(r,\theta,\phi) = \sum_{m = -\infty}^{\infty} f_{m}(r, \theta) e^{i m \phi},
\end{equation}
where each $f_{m}(r, \theta)$ is the $m$-th Fourier mode.  Decomposition into Fourier modes leaves us free to choose the basis to represent the $(r, \theta)$
dependence of arbitrary fields.  In the sphere the most common choice is to use spherical harmonics for the angular $\theta$ dependence \cite{Boyd_Yu_2011}.
These basis functions are a complete set of orthonormal modes on the sphere and behave regularly in the presence of the coordinate singularities
at the poles, $\theta = 0, \pi$.  These singularities take the same form as those in the two-dimensional disk.  In polar coordinates ${ (s, \phi) }$
the $m$-th Fourier mode must behave as
\begin{equation}
f_{m}(s) \sim s^{\left| m \right|} F(s^{2}) \text{ as } s \to 0,
\end{equation}
where $F$ is a well-behaved function of $s^{2}$.  In the sphere $\theta$ plays a role analogous to $s$.
Considering $\sin{\left(\theta\right)} \sim \theta$ and $\cos{\left(\theta\right)} \sim 1-\theta^{2}/2$ as $\theta \to 0$,
the form of a spherical harmonic of degree $l$ exactly matches the regularity requirement above:
\begin{equation}
Y_{l,m}(\theta, \phi) \sim \sin^{\left| m \right|} \left(\theta\right) P_{l,m} \left( \cos{\left(\theta\right)} \right) e^{i m \phi}.
\end{equation}

The final ingredient to standard sphere codes that demand spectral convergence is to expand the radial $r$ dependence in orthogonal polynomials.
The numerical method must appropriately handle the coordinate singularity at the origin, ${ r = 0 }$.
For a spherical harmonic of degree $l$, regularity enforces the condition
\begin{equation}
f_{l,m}(r) \sim r^{l} F(r^{2}) \text{ as } r \to 0.
\end{equation}
Many studies make different choices for the radial dependence \cite{Boyd_Yu_2011, Vasil_Lecoanet_Burns_Oishi_Brown_2019, Matsushima_Marcus_1995}.
Chebyshev polynomials are one important class of radial basis \cite{Boyd_2001, Boyd_Yu_2011, Livermore_Jones_Worland_2007}.
However, without explicit control, singularities may arise in higher derivatives; we contend, as others, that coordinate singularities
are best handled explicitly in the basis functions themselves.  For example, Zernike-type bases play this role
\cite{Vasil_Burns_Lecoanet_Olver_Brown_Oishi_2016, Vasil_Lecoanet_Burns_Oishi_Brown_2019}, with their explicit $r^{l}$ scaling.
The three cascaded transforms - Fourier, spherical harmonics, then Zernike - provide a spectral transformation for
well-behaved scalar fields in the ball.  They are particularly well-suited for fields with rotational invariance - isotropic in $\theta$ and $\phi$.
These bases are sub-optimal for fields that break rotational invariance through the imposition of a preferred ``gyroscopic'' axis of alignment.
For example, the internal structure of stars \cite{Miesch_2005}, giant and minor planets \cite{Stevenson_2003,Schubert_2015}
and off-world oceans \cite{Kivelson_2000,Thomson_Delaney_2001,Gissinger_Petitdemange_2019} are all greatly impacted by rotation about a preferred axis
through the Coriolis force.
Rotationally constrained fluids have a broken symmetry - here the $z$ axis - resulting from a primary (geostrophic) balance between the rotational force and the 
pressure gradient force in the incompressible Navier-Stokes equations.  Flows align cylindrically as we increase rotation rate to meet the
extreme demands of geophysical regimes.
This is a consequence of the Taylor-Proudman theorem, which yields spatially anisotropic, nearly invariant flows in the axial direction \cite{taylor1923experiments,Greenspan_1968}.
Basis functions using spherical harmonics have a rotational invariance not present in these geophysical flows.  We forgo this
traditional approach and seek a new formulation of spectral codes in the sphere.

In this paper we tailor a coordinate system to the cylindrical
structures found in rotating spherical flows and investigate the fully regular orthogonal basis that it inspires.  
The key to the approach is to begin with cylindrical coordinates, then stretch the top and bottom of the cylinder
onto the upper and lower surfaces of the sphere.  This coordinate transformation acts as the backbone for defining
a spectral basis.  The benefit to working in these stretched coordinates is that the coordinate lines match the
structures found in gyroscopically aligned flows.  This implies that, along with a suitable choice of basis, we
can represent geophysical flows with relatively few degrees of freedom compared to spherical harmonics expansions.
By designing a hierarchy of bases we explicitly handle coordinate singularities and implement all calculus operations
required for fluids problems with a sparse system of equations.
Authors have employed this strategy of tuning orthogonal bases to the geometry for other fundamental domains, including on triangles, wedges
and disk slices \cite{Olver_Townsend_Vasil_2019,olver2019orthogonal,olver2020recurrence,Snowball_Olver_2020a}, on quadratic and cubic curves \cite{Olver_Xu_2021,Fasondini_Olver_Xu_2020},
in and on quadratic surfaces of revolution \cite{Olver_Xu_2020} and on spherical caps \cite{Snowball_Olver_2020b}.  In all cases
the geometric volume element induces a hierarchy of Hilbert spaces and orthogonal bases.  Careful choice of domain and codomain of linear operators
maximizes sparsity of the matrix system.  Our work extends this approach to the gyroscopic coordinate system in the full sphere.

In Section~\ref{sec:spherindrical_coordinates} we describe the gyroscopically aligned coordinate system for the sphere and in
Section~\ref{sec:basis} we define a hierarchy of basis functions used to represent scalar and vector fields.  Section~\ref{sec:discretization}
details how the basis naturally leads to sparse matrix operators for all calculus operators needed in fluid dynamics.
We put the basis to the test in Section~\ref{sec:eigenproblems} where we solve three eigenvalue problems from fluids
applications.  We wrap up the paper in Section~\ref{sec:conclusion}.

\section{The Stretched Cylindrical Coordinate System}\label{sec:spherindrical_coordinates}
\begin{figure}[H]
\centering
\includegraphics[width=\linewidth]{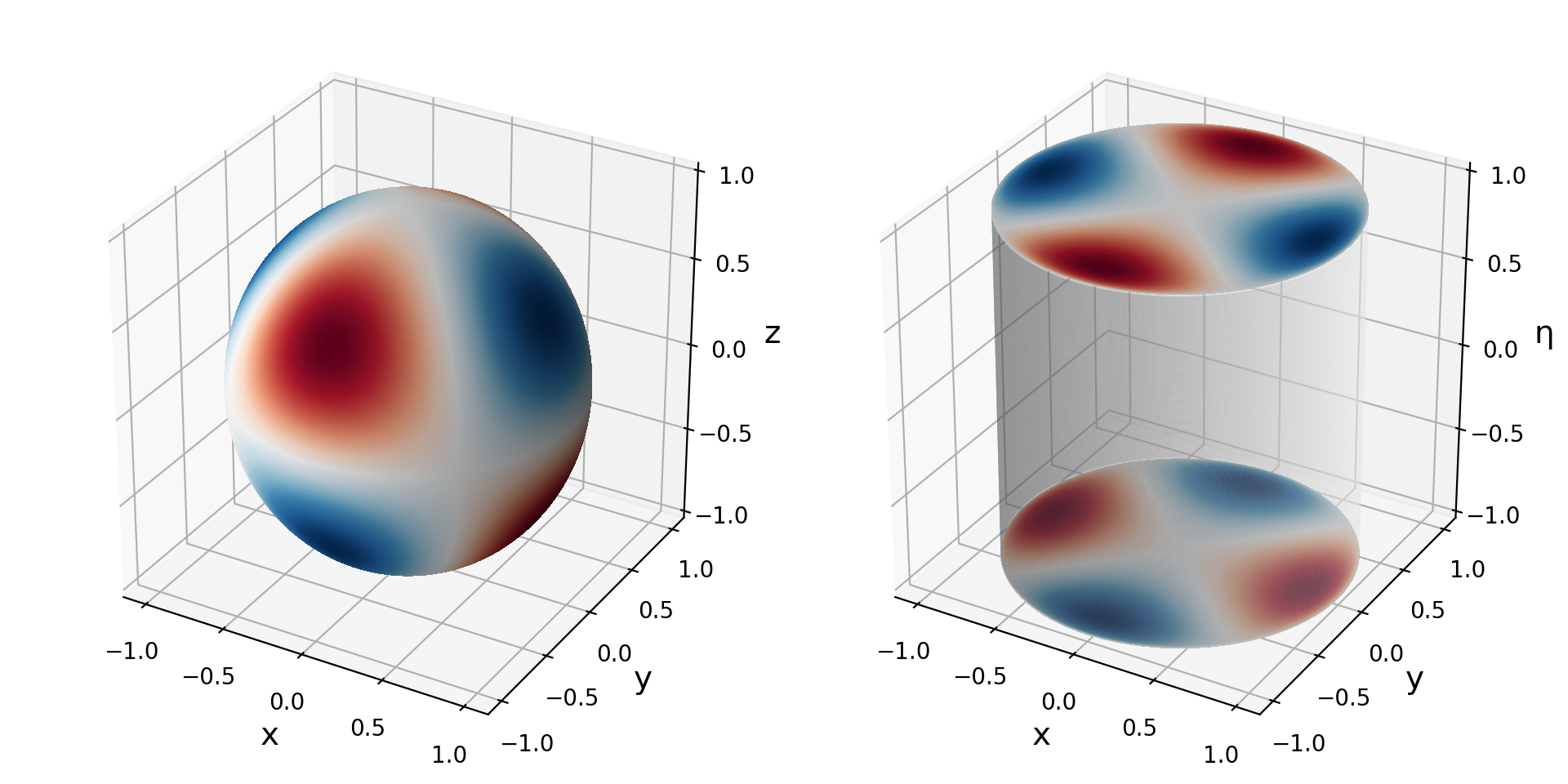}
\caption{The spherical harmonic $Y_{3,2}$ plotted in Cartesian coordinates (left) and in the stretched $\eta$ coordinate (right).
    The upper (respectively, lower) surface of the unit sphere is mapped onto the top (bottom) of the cylinder in the stretched
    coordinate system.  The equator of the sphere is mapped to the entire cylindrical surface at $s = 1$.}
\label{fig:spherinder_mapping}
\end{figure}

We adopt the stretched cylindrical coordinate system ${ (s, \phi, \eta) }$ that maps a cylinder of height two and
unit radius onto a sphere of unit radius, and name the system ``spherindrical coordinates.''
The transformation to Cartesian ${(x, y, z)}$ coordinates is given by
\begin{equation}
\begin{aligned}
x &= s \cos \phi \\
y &= s \sin \phi \\
z &= \eta \sqrt{1 - s^2},
\end{aligned}
\end{equation}
where
\begin{equation}
s \in [0, 1], \hspace{4ex} \phi \in [0, 2 \pi), \hspace{4ex} \eta \in [-1,1].
\end{equation}
Figure~\ref{fig:spherinder_mapping} shows how the surface of the sphere maps onto the upper and lower surfaces of the
cylinder in the stretched coordinates.
The surface equator lives at the cylindrical wall $s = 1$ while the upper (respectively, lower) boundary of the sphere is at ${ \eta = +1 }$ (${ \eta = -1 }$).
Denoting the standard cylindrical coordinates ${ (S, \Phi, Z) }$, the partial derivatives transform as
\begin{equation}
\begin{aligned}
\partial_{S} &= \partial_{s} + \frac{s}{1-s^2} \eta \partial_{\eta} \\
\partial_{\Phi} &= \partial_{\phi} \\
\partial_{Z} &= \frac{1}{\sqrt{1-s^2}} \partial_{\eta}.
\label{eqn:partial_derivative_transformation}
\end{aligned}
\end{equation}
The coupling of the partial derivatives demonstrates that the coordinate vectors $\partial_{s}$ and $\partial_{\eta}$ aren't
orthogonal.  We are thus trading away the decoupled spatial derivatives of spherical harmonics for a sparser representation of
gyroscopically aligned flow morphology in the stretched coordinate system.

\begin{figure}[H]
\centering
\includegraphics[width=.4\linewidth]{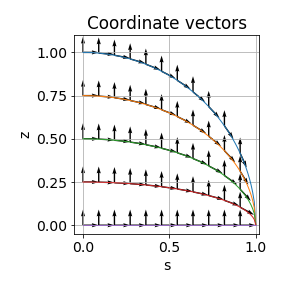}
\includegraphics[width=.4\linewidth]{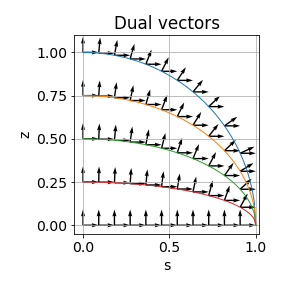}
\caption{Normalized coordinate vectors (left) and their corresponding dual vectors (right).  Each curve is a constant-$\eta$ slice}
\label{fig:coordinate_vectors}
\end{figure}

Figure~\ref{fig:coordinate_vectors} displays the coordinate vectors ${ \partial_s }$ and ${ \partial_\eta }$ and their corresponding dual vectors.
The coordinate singularity at the equator manifests itself in the
convergence of coordinate lines of constant $\eta$ and the linear dependence of the coordinate vectors at $s = 1$.
Due to this behavior of the spherindrical coordinate vectors we elect to represent vector fields in the cylindrical
coordinate basis ${( \vec{\hat{e}}_{S}, \vec{\hat{e}}_{\Phi}, \vec{\hat{e}}_{Z}) }$ but expressed as functions of the spherindrical coordinates ${ (s, \phi, \eta) }$.

\section{The Basis}\label{sec:basis}

\subsection{Basis Functions}
There are two regularity requirements for basis functions in spherindrical coordinates.  The well-known disk singularity implies the $s$ expansion
for azimuthal mode $m$ must take the form
\begin{equation}
f \sim e^{i m \phi} s^{\left| m \right|} F\left(s^2\right) \hspace{4ex} s \to 0,
\end{equation}
where $F$ is any analytic function.  The equatorial singularity induces one further restriction on basis functions - namely we require sufficient
decay as $s \to 1$.  For a vertical monomial in $z$ of degree $l$ we have
\begin{equation}
z^{l} = \left( \eta \sqrt{1-s^{2}} \right)^{l} = \eta^{l} \left( 1-s^{2} \right)^{\frac{l}{2}}
\label{eqn:z-monomial}
\end{equation}
and hence we incorporate the stretching factor $\left( 1-s^{2} \right)^{\frac{l}{2}}$ into our basis functions.  This guarantees our basis is expressible as a
Cartesian polynomial.  See Appendix B for further details on the regularity requirements for a scalar field in the spherinder.
For instance we establish that (\ref{eqn:z-monomial}) is consistent with the regularity conditions for the spherical Laplacian.
Putting the above two constraints together we have, for azimuthal mode $m$ and vertical degree $l$, a smooth function must behave as
\begin{equation}
f(s,\phi,\eta) \sim e^{i m \phi} s^{\left| m \right|} \left(1-s^2\right)^{\frac{l}{2}} P_{l}(\eta) F(s^2),
\end{equation}
where $P_{l}(\eta)$ represents a polynomial of degree $l$ in $\eta$ and $F$ is any analytic function.
Define the coordinate ${ t = 2s^2 - 1 }$ so that ${ -1 \le t \le 1 }$ and
\begin{equation}
f(t,\phi,\eta) \sim e^{i m \phi} (1+t)^{\frac{\left| m \right|}{2}} \left(1-t\right)^{\frac{l}{2}} P_{l}(\eta) \widetilde{F}(t).
\label{eqn:scalar-basis}
\end{equation}
This change of variables ensures our radial dependence $\widetilde{F}$ is even in $s$ and transforms the domain of the radial coordinate
to the natural setting $[-1,1]$ for orthogonal polynomials.
    
The mode in (\ref{eqn:scalar-basis}) suffices for scalar fields in the sphere but cannot capture leading order behavior of vector
fields.  Vector fields have components that decay like ${s^{\left| m \right| \pm 1}}$ as $s \to 0$, which can be seen for example
by taking the gradient of (\ref{eqn:scalar-basis}).  For this reason we introduce the \emph{spin-weight} ${ \sigma \in \{ -1, 0, +1 \} }$,
further developed in Section~\ref{subsec:spinor_basis}, that parameterizes the various decay rates for vector fields along the $z$ axis.

This motivates the definition of a family of orthonormal, spin-weighted basis function parameterized by ${ \alpha \in \mathbb{R} }$, ${ \alpha > -1 }$,
and indexed by integers ${ m \in \mathbb{Z} }$, and ${ l, k \in \mathbb{Z}_{\ge 0} }$.  This \emph{spherinder basis} takes the explicit form
\begin{equation}
\Psi_{m,l,k}^{\sigma,\alpha}(t, \eta, \phi) = \frac{1}{N_{\alpha}} e^{i m \phi} (1 + t)^{\frac{\left| m \right| + \sigma}{2}} (1 - t)^{\frac{l}{2}} P_{l}^{(\alpha, \alpha)}(\eta) P_{k}^{(l+\alpha+\frac{1}{2}, \left| m \right|+\sigma)}(t)
\label{eqn:spherinder_basis}
\end{equation}
where ${ N_{\alpha} = \sqrt{\frac{2 \pi}{2^{2 + \alpha + \frac{1}{2}}}} }$.
Here $P_{n}^{(a,b)}$ is the degree-$n$ Jacobi polynomial orthonormal under the weight function ${ (1-t)^{a} (1+t)^{b} }$.
Each spin-0 basis function $\Psi_{m,l,k}^{0,\alpha}$ is a degree-$\left( m + l + 2k \right)$ homogeneous Cartesian ${(x,y,z)}$ polynomial.
Requiring a homogeneous polynomial representation rules out the more general vertical basis $P_{l}^{(\alpha,\beta)}(\eta)$, which
loses its symmetry about $\eta = 0$ when $\alpha \ne \beta$.

The volume element in spherindrical coordinates is ${ \dd V = s \sqrt{1-s^2} \dd s \dd \phi \dd \eta }$.
We define the family of volume measures,
\begin{equation}
\begin{aligned}
\dd \mu(\alpha) &= \left( 1 - r^2 \right)^{\alpha} \dd V \\
                &= \left( 1 - \eta^2 \right)^{\alpha} \left( 1 - s^2 \right)^{\alpha+\frac{1}{2}} s \dd s \dd \phi \dd \eta \\
                &= \frac{1}{2^{2 + \alpha + \frac{1}{2}}} \left( 1 - \eta^2 \right)^{\alpha} \left( 1 - t \right)^{\alpha+\frac{1}{2}} \dd t \dd \phi \dd\eta,
\end{aligned}
\end{equation}
to produce a hierarchy of bases in the parameter $\alpha$.
These basis functions are orthonormal polynomials under the inner product
\begin{equation}
\begin{aligned}
\left \langle \Psi_{m,l,k}^{\sigma, \alpha} \hspace{1ex}, \Psi_{m',l',k'}^{\sigma, \alpha} \right \rangle_{(\alpha)} 
    &= \int \Psi_{m,l,k}^{\sigma, \alpha} \Psi_{m',l',k'}^{\sigma, \alpha} \dd \mu(\alpha) \\
    &= \frac{1}{2^{2 + \alpha + \frac{1}{2}}} \int \Psi_{m,l,k}^{\sigma, \alpha} \Psi_{m',l',k'}^{\sigma, \alpha} (1-\eta^2)^{\alpha} \left(1-t\right)^{\alpha + \frac{1}{2}} \dd t \dd \phi \dd \eta \\
    &= \frac{1}{2 \pi} \int_{0}^{2 \pi} e^{-i (m-m')\phi} \dd \phi \int_{-1}^{1} P_{l}^{(\alpha,\alpha)}(\eta) P_{l'}^{(\alpha,\alpha)}(\eta) (1-\eta^2)^{\alpha} \dd \eta \\
&\hspace{4ex} \times             \int_{-1}^{1} P_{k}^{(l+\alpha+\frac{1}{2},\left|m\right|+\sigma)}(t) P_{k'}^{(l+\alpha+\frac{1}{2},\left|m\right|+\sigma)}(t) \left( 1-t \right)^{l + \alpha + \frac{1}{2}} \left( 1+t \right)^{\left|m\right|+\sigma} \dd t \\
    &= \delta_{m,m'} \delta_{l,l'} \delta_{k,k'},
\label{eqn:inner_product}
\end{aligned}
\end{equation}
where $\delta_{a,b}$ is the Kronecker delta.

We render Jacobi polynomial calculus sparse by allowing the Jacobi parameters to float naturally under the action of operators.
The derivative of a Chebyshev polynomial of the first kind, $T_{k}$, represented in the $T_{n}$ basis is an upper triangular operator.
Representing the derivative instead in the basis of Chebyshev polynomials of the second kind, $U_{n}$, diagonalizes the operator. 
This corresponds to raising both Jacobi parameters - from ${\left( -1/2, -1/2 \right)}$ to ${\left( 1/2, 1/2 \right)}$ - under the action of the derivative.
This key observation allows us to construct a sparse calculus for PDEs in the sphere; we introduce the parameter ${ \alpha > -1 }$ into our basis functions and,
along with the parameter $\sigma$, we have a sufficiently general hierarchy of basis functions to represent regular scalar and
vector fields in the ball.

We note that though derivatives increment the $\alpha$ index to maintain sparsity of the matrix system,
it is rather arbitrary.  To match the volume element in the $\eta$
coordinate we typically choose ${ \alpha = 0 }$ as the starting point for our discretization.  This corresponds to Legendre
polynomials in the vertical direction with a unit integration weight.  Another useful candidate is Chebyshev $T_{n}$
polynomials for the vertical expansion, for which we set ${ \alpha = -\frac{1}{2} }$.

\subsection{Zoology of Basis Functions}
\begin{figure}[H]
\centering
\includegraphics[width=\linewidth]{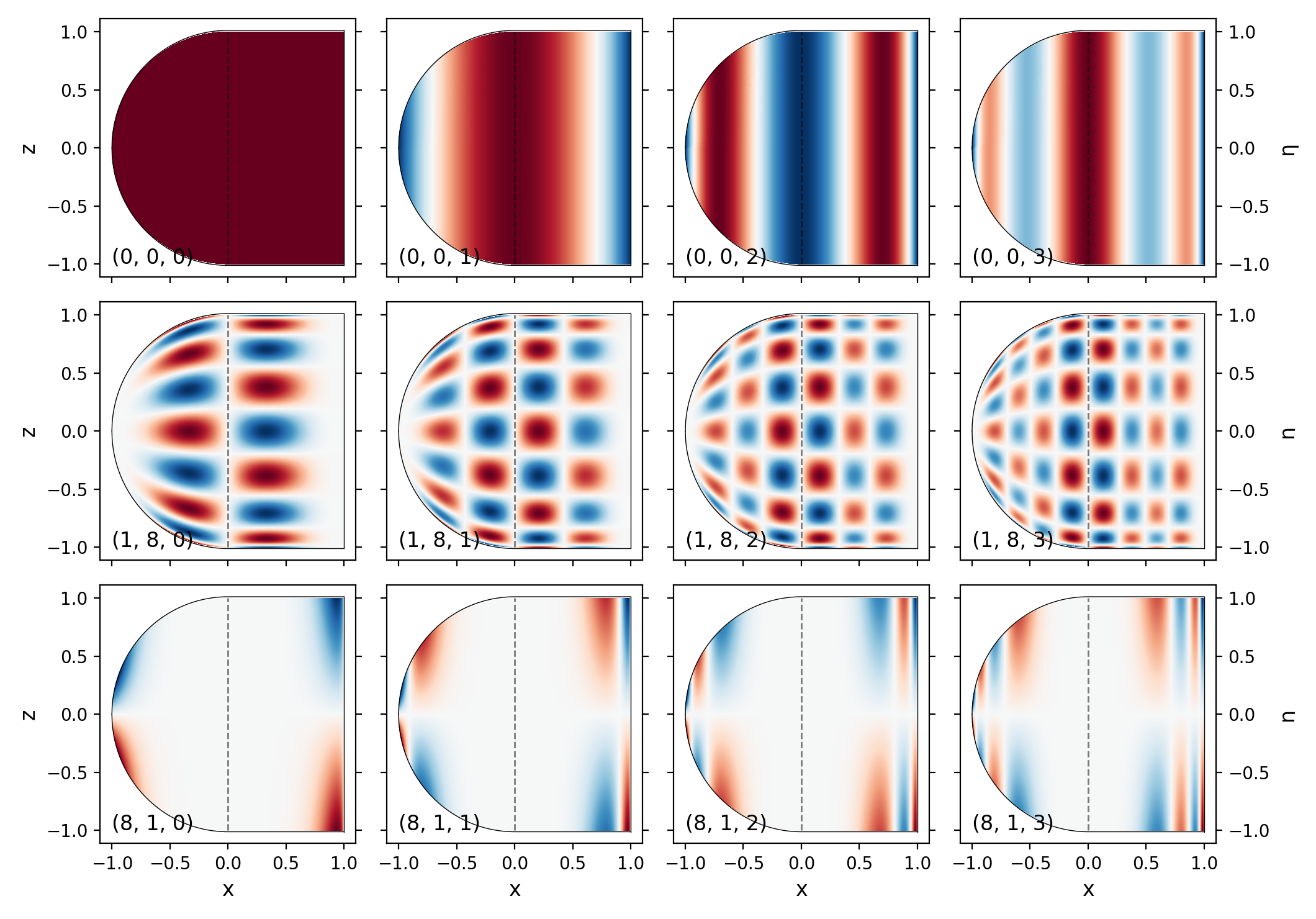}
\caption{Meridional slice of basis functions $\Psi_{m,l,k}^{0,-\frac{1}{2}}$ at $\phi = 0$ (see (\ref{eqn:spherinder_basis})).  
    The left half ($x < 0$) of each subplot uses the $z$ coordinate while the right half ($x > 0$)
    uses the $\eta$ coordinate.  This provides views of the basis functions both in their natural coordinates and in their physical geometry.
    The $(m,l,k)$ index for each basis function appears at the bottom left of each plot.
    Each column represents a different radial mode.  The rows are
    distinct choices of the $m$ and $l$ indices.  Notice in the high $m$ basis functions (bottom row) the $s^{\left| m \right|}$
    decay at the $z$ axis (vertical center line).  Visual decay as $s \to 1$ emerges in the high $l$ basis functions (middle row).  Basis functions
    shown are all ${ \sigma = 0 }$, but rather than show the natural ${ \alpha = 0 }$ bases we choose ${ \alpha = -1/2 }$.
    These basis functions are Chebyshev $T_{n}$ polynomials in $\eta$ and the more uniform oscillations show better contrast when plotting.}
\label{fig:basis-zoo}
\end{figure}

Figure~\ref{fig:basis-zoo} exhibits some of these basis functions at ${ \phi = 0 }$ for various vertical degrees $l$ and radial degrees $k$.
A key feature of the basis functions is their behavior at the coordinate singularities - the $z$ axis ($s=0$) and the equator ($s=1$).
All scalar fields decay like $s^{\left| m \right|}$ as $s \to 0$ for azimuthal degree $m$ and like ${ (1 - s^2)^{\frac{l}{2}} }$ as $s \to 1$ for vertical degree $l$.
We denote the polynomial degree in the vertical direction by $l$ and the degree in the radial direction by $k$.

\subsection{Spinor Basis}\label{subsec:spinor_basis}
The coordinate vector basis ${ (\partial_{s}, \partial_{\phi}, \partial_{\eta}) }$ is not directly used for interpreting
physical problems.  For instance, at the equator these vectors become linearly dependent, a manifestation of the coordinate
singularity there.  We instead utilize the unit-normalized cylindrical coordinate vectors ${ (\vec{\hat{e}}_{S}, \vec{\hat{e}}_{\Phi}, \vec{\hat{e}}_{Z}) }$
to define an orthonormal basis for representing vector fields. 

Throughout the paper we will denote the \emph{spin weight} of a component of a vector field by ${ \sigma \in \{ -1, 0, +1 \} }$.
To decouple the coordinate vectors under gradient action we, following \cite{Vasil_Burns_Lecoanet_Olver_Brown_Oishi_2016,Vasil_Lecoanet_Burns_Oishi_Brown_2019}, define the spinor basis ${ \left( \vec{\hat{e}}_{\sigma} \right) }$:
\begin{equation}
\vec{\hat{e}}_{\pm} = \frac{1}{\sqrt{2}} \left(\vec{\hat{e}}_{S} \mp i \vec{\hat{e}}_{\Phi} \right)
\hspace{4ex} \implies \vec{\hat{e}}_{S} = \frac{1}{\sqrt{2}} \left( \vec{\hat{e}}_{+} + \vec{\hat{e}}_{-} \right),
\hspace{4ex} \vec{\hat{e}}_{\Phi} = \frac{i}{\sqrt{2}} \left( \vec{\hat{e}}_{+} - \vec{\hat{e}}_{-} \right),
\label{eqn:spinor_basis}
\end{equation}
with $\vec{\hat{e}}_{0} \triangleq \vec{\hat{e}}_{Z}$.  Then we represent a vector field $\vec{u}(s,\phi,\eta)$ as
\begin{equation}
\vec{u}(s,\phi,\eta) = u \vec{\hat{e}}_{S} + v \vec{\hat{e}}_{\Phi} + w \vec{\hat{e}}_{Z} = u_{+} \vec{\hat{e}}_{+} + u_{-} \vec{\hat{e}}_{-} + u_{0} \vec{\hat{e}}_0, 
\end{equation}
where $u$, $v$, $w$ and $u_{\sigma}$ are all functions of the stretched coordinates $(s, \phi, \eta)$.
We find
\begin{equation}
u_{\pm} = \frac{1}{\sqrt{2}} \left( u \pm i v \right) \\
\hspace{4ex} \implies u = \frac{1}{\sqrt{2}} \left( u_{+} + u_{-} \right),
\hspace{4ex} v = -\frac{i}{\sqrt{2}} \left( u_{+} - u_{-} \right)
\end{equation}
and $u_{0} \equiv w$.  We then have
\begin{equation}
\vec{\hat{e}}_{+} \cdot \vec{\hat{e}}_{+} = \vec{\hat{e}}_{-} \cdot \vec{\hat{e}}_{-} = 0, 
\hspace{4ex} \vec{\hat{e}}_{+} \cdot \vec{\hat{e}}_{-} = \vec{\hat{e}}_{-} \cdot \vec{\hat{e}}_{+} = 1.
\label{eqn:spinor-dots}
\end{equation}
The results (\ref{eqn:spinor-dots}) imply the dual basis to the spinor basis is given by the complex conjugate.
Hence to extract a spin component $u_{\sigma}$ from a vector $\vec{u}$ we have
\begin{equation}
u_{\sigma} = \left( \vec{\hat{e}}_{\sigma} \right)^{\dagger} \cdot \vec{u} = \vec{\hat{e}}_{\sigma}^{*} \cdot \vec{u} = \vec{\hat{e}}_{-\sigma} \cdot \vec{u},
\end{equation}
where the dagger denotes the dual vector and the star denotes complex conjugation.

The cross product with ${ \vec{\hat{e}}_{Z} \equiv \vec{\hat{e}}_{0} }$ then gives
\begin{equation}
\vec{\hat{e}}_{0} \times \vec{\hat{e}}_{\pm} = \pm i \vec{\hat{e}}_{\pm}.
\end{equation}

The horizontal gradient takes the form 
\begin{equation}
\nabla_{\perp} = \vec{\hat{e}}_{S} \nabla_{S} + \vec{\hat{e}}_{\Phi} \nabla_{\Phi} = \vec{\hat{e}}_{-} \nabla_{-} + \vec{\hat{e}}_{+} \nabla_{+},
\end{equation}
where $\nabla_{\pm}$ acts on a single azimuthal mode ${ e^{i m \phi} f_m(s, \eta) }$ by
\begin{equation}
\nabla_{\pm} \equiv \left[ \partial_{S} \mp \frac{m}{S} \right].
\end{equation}
Then the gradient connection is diagonal:
\begin{equation}
\nabla_{\pm} \vec{\hat{e}}_{+} = \mp \frac{1}{S} \vec{\hat{e}}_{+}, \hspace{4ex}
\nabla_{\pm} \vec{\hat{e}}_{-} = \pm \frac{1}{S} \vec{\hat{e}}_{-}.
\end{equation}
Critically, this diagonalization means the vector Laplacian operator doesn't couple the vector components
represented in the spinor basis, thereby improving the sparsity of calculus operations.

Another important motivation for the spinor basis is that the spin components behave predictably at the $z$ axis.
For ${ \sigma \in \{ -1, 0, +1 \} }$, the $u_{\sigma}$ component of the vector field decays like
\begin{equation}
u_{\sigma} = \vec{\hat{e}}_{\sigma}^{*} \cdot \vec{u} \sim s^{\left| m \right| + \sigma } F_{\sigma} \left( s^{2} \right), \hspace{4ex} s \to 0,
\end{equation}
where $F_{\sigma}$ is an arbitrary well-behaved function of $s^{2}$.  This means the the basis functions $\Psi_{m,l,k}^{\sigma,\alpha}$ defined
in (\ref{eqn:spherinder_basis}) behave precisely as needed to represent vector fields regular throughout the ball.

\subsection{Field Expansions}
With Equation~(\ref{eqn:spherinder_basis}) and the results of Section~\ref{subsec:spinor_basis}
we have the definitions required to represent scalars and vector fields.  We decompose a scalar field as
\begin{equation}
f(t = 2 s^2 - 1, \phi, \eta) = \frac{1}{N_{\alpha}} \sum_{m = -\infty}^{\infty} e^{i m \phi} (1+t)^{\frac{\left| m \right|}{2}} 
    \sum_{l = 0}^{\infty} (1-t)^{\frac{l}{2}} P_{l}^{(\alpha,\alpha)}(\eta) 
        \sum_{k = 0}^{\infty} P_{k}^{(l+\alpha+\frac{1}{2},\left| m \right|)}(t) \widehat{F}_{m,l,k},
\label{eqn:scalar-expansion}
\end{equation}
where $\alpha$ can be freely chosen.  We represent vector fields as a sum of their spin components using
\begin{equation}
\vec{u}(t,\phi,\eta) = \frac{1}{N_{\alpha}} \sum_{\sigma} \vec{\hat{e}}_{\sigma} \sum_{m = -\infty}^{\infty} e^{i m \phi} (1+t)^{\frac{\left| m \right|+\sigma}{2}} 
    \sum_{l = 0}^{\infty} (1-t)^{\frac{l}{2}} P_{l}^{(\alpha,\alpha)}(\eta) 
        \sum_{k = 0}^{\infty} P_{k}^{(l+\alpha+\frac{1}{2},\left| m \right|+\sigma)}(t) \widehat{U}_{m,l,k}^{\sigma}.
\end{equation}
Using Jacobi polynomial algebra we obtain sparse calculus operations on fields expanded in the spherindrical basis.  Appendix A details the necessary operators
for most PDEs, for example those occurring in fluid dynamics.  The key takeaway is that, though operations typically couple vertical and radial modes ($l$ and $k$ indices), the expansion
enables a sparse representation of gyroscopically aligned flows while explicitly conforming to coordinate singularities.

In what follows we drop the ${m,l,k}$ subscript and use $\widehat{F}$ and $\widehat{U}^{\sigma}$ to denote the set of all spectral coefficients for scalar
and vector fields, respectively.

\section{Discretization}\label{sec:discretization}
We now demonstrate the numerical constructions required to use the basis to solve PDE systems in a spherical ball.  We first
give an overview of how calculus operators act on the basis functions, then demonstrate how to employ these in the solution
to a simple test problem.

\subsection{Operators}\label{sec:operators}
Jacobi polynomial algebra allows us to implement calculus operators on our basis functions with sparse matrix operations.
For clarity we relegate computation of the operator coefficients to Appendix A; here we simply demonstrate the results.  Operators typically map not only between
the ${ (l,k) }$ indices of the basis functions, but also may move ${ (\sigma, \alpha) }$ according to the operation.  We
begin by listing the action of differential operators, then demonstrate the operators involving the spherical radial coordinate $r$.
For what follows we make use of the inner product (\ref{eqn:inner_product}) to define a Hilbert space norm
\begin{equation}
\norm{ f }_{(\alpha)}^2 = \left\langle f, f \right\rangle_{(\alpha)}.
\label{eqn:norm-def}
\end{equation}

When we define operators in the subsequent sections we will use the
superscript $\delta$ to denote the \emph{change in spin-weight} $\sigma$ due to applying the operator.  We emphasize that scalar fields live in spin-weight
${ \sigma = 0 }$ while vectors naturally decompose into spins ${ \sigma \in \{ -1, 0, +1 \} }$.  Operators that map scalars to vectors
will split into three parts that take spin 0 to the appropriate spin weights.  Likewise operators that map vectors to scalars will split
into a $\sigma$-lowering operator (${ \delta = -1 }$) acting on the ${ \sigma = +1 }$ vector component, a $\sigma$-raising operator (${ \delta = +1 }$) acting on the ${ \sigma = -1 }$ vector
component and a $\sigma$-preserving operator (${ \delta = 0 }$) acting on the ${ \sigma = 0 }$ component.

\subsubsection{Regularity}
As previously discussed, regularity in spherindrical coordinates is associated with the coordinate $s$ at ${s = 0}$ and ${s = 1}$.
To capture the regularity structure of fields in the ball we define the \emph{regularity space} of degree $m$ as
\begin{equation}
\textrm{Reg}(m) = \left\{ f:[0,1] \to \mathbb{C} \hspace{2ex} \text{ s.t. } \hspace{2ex} 
    f(s) \sim s^{m} F\left(s^2\right) \hspace{1ex} \text{ as } \hspace{1ex} s \to 0
    \right\},
\end{equation}
where $F(s^2)$ is any even function of $s$ that is analytic in neighborhood of $s = 0$.
Then for the $m$-th azimuthal mode scalar fields live in $\textrm{Reg}(\left| m \right|)$ while vector fields decompose into 
the direct sum of three regularity spaces, $\textrm{Reg}(\left| m \right|-1) \oplus \textrm{Reg}(\left| m \right|) \oplus \textrm{Reg}(\left| m \right| +1)$
\cite{Vasil_Burns_Lecoanet_Olver_Brown_Oishi_2016}, corresponding to the spin components of the vector field.

Scalars and vector fields (as well as higher rank tensor fields)
expressed in the spin-cylindrical basis $(\vec{\hat{e}}_{\sigma})$ always decay like $\left(1-s^2\right)^{\frac{l}{2}}$ as $s \to 1$ for a monomial in $\eta$ of
degree $l$.  Noting that the Cartesian polynomials are a complete basis for functions in the sphere and, using the relation ${ z^{l} = \eta^{l} \left( 1 - s^{2} \right)^{\frac{l}{2}} }$,
any function ${ g(s,\eta) : [0,1] \times [-1,1] }$ regular of degree $m$ in its first argument can be decomposed as
\begin{equation}
g(s,\eta) = \sum_{l=0}^{\infty}{ \eta^{l} \left(1 - s^{2} \right)^{\frac{l}{2}} g_{l}(s) }
\end{equation}
for functions $g_{l} \in \textrm{Reg}(m)$.  The decomposition above matches the form of our modal expansion for scalar fields (\ref{eqn:scalar-expansion}),
though in this form we use the monomial basis for simplicity of exposition rather than a basis of orthogonal polynomials in the $\eta$ coordinate.

We now define a hierarchy of Hilbert spaces indexed by real parameter $\alpha > -1$:
\begin{equation}
\mathcal{H}^{\alpha}(m) = \left\{ f \in \textrm{Reg}(m) \hspace{2ex}\text{ s.t. }\hspace{2ex} \norm{f}_{(\alpha)} < \infty \right\},
\end{equation}
where $\norm{.}_{(\alpha)}$ refers to the norm induced by the inner product (\ref{eqn:inner_product}).
Differential operators map between these Hilbert spaces.  Jacobi polynomial algebra provides a sparse
representation of these operators acting on the spherinder basis.

\subsubsection{Differential Operators}
Recall that ${ \alpha > -1 }$ is a numerical index.  This parameter defines the Jacobi polynomial class
for the vertical ($\eta$) and radial (${ t = 2 s^{2} - 1 }$) parts of the basis functions.  Choosing $\alpha = 0$ utilizes the
geometric volume element and results in Legendre polynomials in the vertical coordinate.  Other choices of
$\alpha$ are possible and will change the type of polynomials used in the basis functions.
Also recall that $\sigma \in \{ -1, 0, +1 \}$ is the spin weight and
${ \delta \in \{ -1, 0, +1 \} }$ is the spin weight increment given by application of an operator.

The \emph{scalar gradient} operator acts on the expansion via
\begin{equation}
\begin{aligned}
\mathcal{G}^{\delta} &: \mathcal{H}^{\alpha}(m) \to \mathcal{H}^{\alpha+1}(m+\delta), \\
\del f &\longleftrightarrow \vec{\hat{e}}_{+} \mathcal{G}^{+} \widehat{F} + \vec{\hat{e}}_{-} \mathcal{G}^{-} \widehat{F} + \vec{\hat{e}}_{0} \mathcal{G}^{0} \widehat{F}.
\end{aligned}
\end{equation}
Note that the gradient naturally decomposes into three spin weights.

\begin{figure}[H]
\centering
\begin{tikzpicture}
    \node (A) at (0,0) {$f \in \mathcal{H}^{\alpha}(m)$};
    \node (B) at (5,+2) {$\vec{\hat{e}}_{+}^{*} \cdot \del f \in \mathcal{H}^{\alpha+1}(m+1) $};
    \node (C) at (5, 0) {$\vec{\hat{e}}_{0}^{*} \cdot \del f \in \mathcal{H}^{\alpha+1}(m) $};
    \node (D) at (5,-2) {$\vec{\hat{e}}_{-}^{*} \cdot \del f \in \mathcal{H}^{\alpha+1}(m-1) $};
    \node[circle, draw=black, scale=0.7] (E) at (10,0) {+};
    \node (F) at (12.5,0) {$\del^{2} f \in \mathcal{H}^{\alpha+2}(m) $};
    \path[->] (A) edge[bend left=10]  node[above] {$\mathcal{G}^{+}$} (B);
    \path[->] (A) edge                node[above] {$\mathcal{G}^{0}$} (C);
    \path[->] (A) edge[bend right=10] node[above] {$\mathcal{G}^{-}$} (D);
    \path[->] (B) edge[bend left=10]  node[above] {$\mathcal{D}^{-}$} (E);
    \path[->] (C) edge                node[above] {$\mathcal{D}^{0}$} (E);
    \path[->] (D) edge[bend right=10] node[above] {$\mathcal{D}^{+}$} (E);
    \path[->] (E) edge node {} (F);
\end{tikzpicture}
\caption{Action of the $\mathcal{G}^{\delta}$ and $\mathcal{D}^{\delta}$ operators on a scalar field $f$
    with respect to regularity spaces.  The scalar Laplacian is computed by cascading these operators
    with the relation ${ \del^{2} f = \del \cdot \left( \del f \right) }$.  Observe how the $\alpha$
    index is incremented each time a derivative is taken, while the scalar gradient decomposes into spin
    weights $\sigma$ which are then recombined through the divergence operator to produce a spin-0 field.}
\label{fig:laplacian-schematic}
\end{figure}
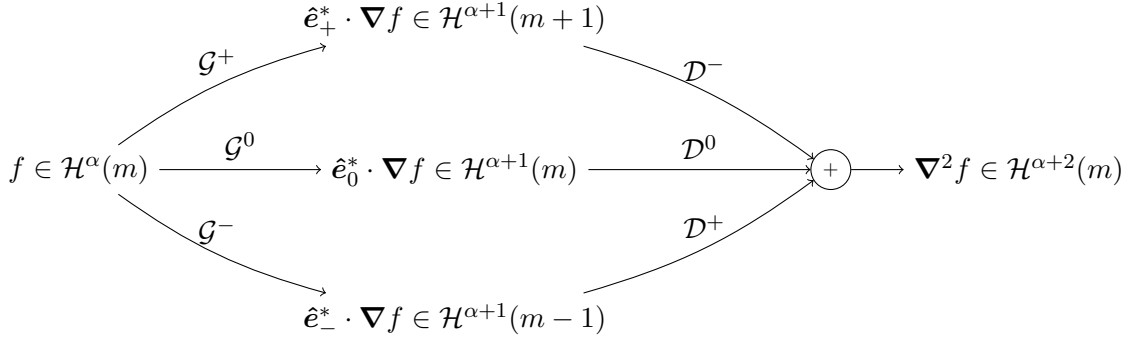

We define the \emph{vector divergence} operator on a vector by
\begin{equation}
\begin{aligned}
\mathcal{D}^{\delta} &: \mathcal{H}^{\alpha}(m-\delta) \to \mathcal{H}^{\alpha+1}(m), \\
\del \cdot \vec{u} &\longleftrightarrow \mathcal{D}^{-} \widehat{U}^{+} + \mathcal{D}^{+} \widehat{U}^{-} + \mathcal{D}^{0} \widehat{U}^{0}.
\end{aligned}
\end{equation}
We consistently define superscripts on the operators to indicate the direction of motion of $\sigma$ - 
namely $\mathcal{D}^{+}$ maps ${ \sigma \mapsto \sigma+1 }$ whereas $\mathcal{D}^{-}$ maps ${ \sigma \mapsto \sigma-1 }$.  
Figure~\ref{fig:laplacian-schematic} depicts a schematic how the gradient and divergence operators work together by first
splitting a scalar into three spin components then recombining them.
We emphasize that these operators map the three spin-weighted Hilbert spaces ${ \mathcal{H}(\{-,0,+\},\alpha) }$ to 
the \emph{same} Hilbert space ${ \mathcal{H}(0,\alpha+1) }$, hence the summation is well-formed.  
This will always be the case - calculus operations on the basis elements are well-defined.  You
can always anticipate finding a sparse representation in terms of Jacobi polynomial algebra.

The \emph{scalar Laplacian} is defined as
\begin{equation}
\begin{aligned}
\mathcal{L} &: \mathcal{H}^{\alpha}(m) \to \mathcal{H}^{\alpha+2}(m), \\
\del^2 f &\longleftrightarrow \mathcal{L} \widehat{F}.
\end{aligned}
\end{equation}
Note the spin weight is unchanged but $\alpha$ is incremented twice, corresponding to the second-order spatial derivatives present in the Laplacian.
We compute the operator from the identity
\begin{equation}
\del^2 f = \del \cdot ( \del f ) \implies \mathcal{L} = \mathcal{D}^{-} \mathcal{G}^{+} + \mathcal{D}^{+} \mathcal{G}^{-} + \mathcal{D}^{0} \mathcal{G}^{0}
\end{equation}
rather than by explicit matrix element computation for non-compound operators in Appendix A.

We define the \emph{vector curl} through the matrix operator
\begin{equation}
\begin{aligned}
\mathcal{C}_{\sigma}^{\delta} &: \mathcal{H}^{\alpha}(m+\sigma-\delta) \to \mathcal{H}^{\alpha+1}(m+\sigma), \\
\del \vec{\times} \vec{u} 
   &\longleftrightarrow \begin{bmatrix} \vec{\hat{e}}_{+} & \vec{\hat{e}}_{-} & \vec{\hat{e}}_{0} \end{bmatrix}
      \begin{bmatrix}
        \mathcal{C}_{+}^{0} & 0 & \mathcal{C}_{+}^{+} \\
        0 & \mathcal{C}_{-}^{0} & \mathcal{C}_{-}^{-} \\
        \mathcal{C}_{0}^{-} & \mathcal{C}_{0}^{+} & 0
      \end{bmatrix} 
      \begin{bmatrix} \widehat{U}^{+} \\ \widehat{U}^{-} \\ \widehat{U}^{0} \end{bmatrix}.
\end{aligned}
\end{equation}

We define the \emph{vector Laplacian} $\del^2$ by the relation
\begin{equation}
\begin{aligned}
\mathcal{L}_{\sigma} &: \mathcal{H}^{\alpha}(m+\sigma) \to \mathcal{H}^{\alpha+2}(m+\sigma), \\
\del^2 \vec{u} &\longleftrightarrow \sum_{\sigma} \vec{\hat{e}}_{\sigma} \mathcal{L}_{\sigma} \widehat{U}^{\sigma}.
\end{aligned}
\end{equation}
We placed the $\sigma$ in the subscript here to emphasize that the operator does not modify $\sigma$ - in fact it is diagonal in the three components
as suggested by the coordinate system.  As for the scalar Laplacian we compute the vector Laplacian from existing operators with the identity
${ \del^2 \vec{u} = \del (\del \cdot \vec{u}) - \del \vec{\times} \del \vec{\times} \vec{u} }$.

\begin{figure}[H]
\centering
\includegraphics[width=0.8\linewidth]{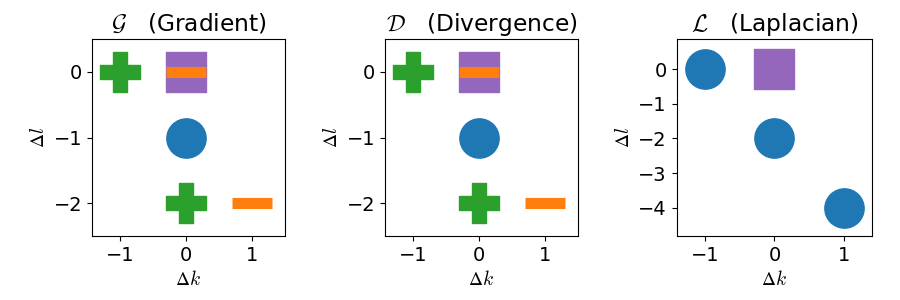}
\caption{Sparse mode coupling for the differential operators.  The purple square represents the input mode.  The remaining symbols denote the outputs
    of the operator.  Green `+' symbols are the outputs of the $\mathcal{G}^{+}$ operator that couple to the $u_{+}$ velocity component.  
    Orange `-' symbols correspond to $\mathcal{G}^{-}$ and $u_{-}$, while the blue circle represents $\mathcal{G}^{0}$ and $u_{0}$.  The input
    mode is assumed to have arbitrary $l,k$ index bounded away from zero.  The vertical axis in the plots represents change in $l$ index,
    and the horizontal axis represents change in $k$ index.  A diagonal operator would have a single output symbol in the $(\Delta l,\Delta k) = (0,0)$ position.
    From the diagram we see differential operators couple not only the $k$ index, but also the $l$ index.  This results from the coupled partial derivatives
    of the non-orthogonal spherindrical coordinate system.}
\label{fig:sparsity_differential}
\end{figure}

Figure~\ref{fig:sparsity_differential} demonstrates mode coupling for each of the operators for a single input mode, denoted at the ${(l,k) = (0,0)}$ position
with a purple square.  In general operators couple up and down in both the $l$ and $k$ indices.  The plot markers portray the change in spin index as defined for 
$\mathcal{G}^{\delta}$ and $\mathcal{D}^{\delta}$ above.  The plus marker corresponds ${\delta = +1}$, the minus marker corresponds to ${\delta = -1}$ and the disk
corresponds to ${\delta = 0}$.

Figure~\ref{fig:sparsity_differential} contains a crucial implication.  Notice the orange minus markers corresponding to the $\mathcal{G}^{-}$
and $\mathcal{D}^{-}$ operators.  In both cases one of the markers \emph{lowers the $l$ index by two} while simultaneously \emph{raising the $k$
index by one}.  If we naively truncate an expansion with the same radial degree $N_{\text{max}}$ for each vertical mode $l$ then we clearly see
$\mathcal{G}^{-} \left(l,N_{\text{max}} \right) \mapsto \left(l-2,N_{\text{max}}+1 \right)$, a basis function outside our basis set!  This means differential
operators cannot be exact if we choose this discretization strategy.  Fortunately there is a simple fix.  Any time a differential operator increases the $k$ index
it occurs in tandem with a decrease in $l$ by at least two.  We therefore choose the maximum radial degree to be a function of $l$:
\begin{equation}
N(l) = N_{\text{max}} - \floor{ \frac{l}{2} },
\end{equation}
where $\floor{ . }$ denotes the floor operator.
This is analogous to the triangular truncation of spherical harmonics and can be interpreted similarly - we must maintain a constant
maximum total polynomial degree for each vertical mode.  Since our radial variable ${t = 2 s^{2} - 1}$ has degree two, the maximum radial degree
must decrease every time the vertical degree increases by two.  In the numerical examples that follow we will always make use of this triangular truncation
strategy.

\subsubsection{Spherical Radial Operators}

Multiplication of a scalar by ${ \vec{r} = r \vec{\hat{e}}_r }$ is given by
\begin{equation}
\begin{aligned}
\mathcal{R}^{\delta} &: \mathcal{H}^{\alpha}(m) \to \mathcal{H}^{\alpha}(m+\delta), \\
r f \vec{\hat{e}}_r &\longleftrightarrow \vec{\hat{e}}_{+} \mathcal{R}^{+} \widehat{F} + \vec{\hat{e}}_{-} \mathcal{R}^{-} \widehat{F} + \vec{\hat{e}}_{0} \mathcal{R}^{0} \widehat{F}.
\end{aligned}
\end{equation}

\begin{figure}[H]
\centering
\begin{tikzpicture}
    \node (A) at (0,0) {$f \in \mathcal{H}^{\alpha}(m)$};
    \node (B) at (5,+2) {$\vec{\hat{e}}_{+}^{*} \cdot \vec{r} f \in \mathcal{H}^{\alpha}(m+1) $};
    \node (C) at (5, 0) {$\vec{\hat{e}}_{0}^{*} \cdot \vec{r} f \in \mathcal{H}^{\alpha}(m) $};
    \node (D) at (5,-2) {$\vec{\hat{e}}_{-}^{*} \cdot \vec{r} f \in \mathcal{H}^{\alpha}(m-1) $};
    \node[circle, draw=black, scale=0.7] (E) at (10,0) {+};
    \node (F) at (12.1,0) {$r^{2} f \in \mathcal{H}^{\alpha}(m) $};
    \path[->] (A) edge[bend left=10]  node[above] {$\mathcal{R}^{+}$} (B);
    \path[->] (A) edge                node[above] {$\mathcal{R}^{0}$} (C);
    \path[->] (A) edge[bend right=10] node[above] {$\mathcal{R}^{-}$} (D);
    \path[->] (B) edge[bend left=10]  node[above] {$\mathcal{E}^{-}$} (E);
    \path[->] (C) edge                node[above] {$\mathcal{E}^{0}$} (E);
    \path[->] (D) edge[bend right=10] node[above] {$\mathcal{E}^{+}$} (E);
    \path[->] (E) edge                node[above] {} (F);
\end{tikzpicture}
\caption{Action of the $\mathcal{R}^{\delta}$ and $\mathcal{E}^{\delta}$ operators on a scalar field $f$
    with respect to regularity spaces.   Observe how the $\alpha$ index is never incremented since no derivatives are taken. 
    This is in contrast to the spin weights $\sigma$ for the vector field $\vec{r} f$ which are recombined into a spin-0 scalar field.}
\label{fig:radial-schematic}
\end{figure}
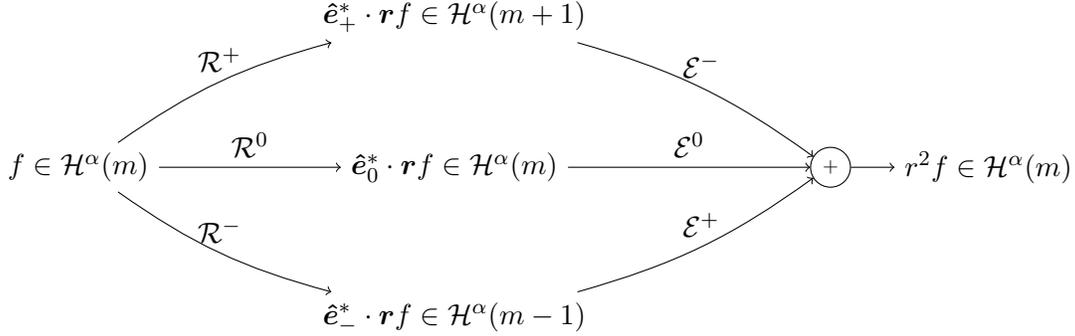

To extract the spherical radial component of a vector field we have
\begin{equation}
\begin{aligned}
\mathcal{E}^{\delta} &: \mathcal{H}^{\alpha}(m-\delta) \to \mathcal{H}^{\alpha}(m), \\
\vec{r} \cdot \vec{u} &\longleftrightarrow \mathcal{E}^{-} \widehat{U}^{+} + \mathcal{E}^{+} \widehat{U}^{-} + \mathcal{E}^{0} \widehat{U}^{0}.
\end{aligned}
\end{equation}
Note that we do not have an operation for ${ \vec{\hat{e}}_r \cdot \vec{u} }$.  This results in division by $r$ which doesn't decouple into
a product of operators in the $t$ and $\eta$ coordinates separately.  Figure~\ref{fig:radial-schematic} shows how the $\mathcal{R}^{\delta}$
and $\mathcal{E}^{\delta}$ operators act on regularity spaces to form the composite operator, multiplication by $r^{2}$.

We define multiplication by ${ 1 - r^2 }$ with the operator
\begin{equation}
\begin{aligned}
\mathcal{S} &: \mathcal{H}^{\alpha}(m) \to \mathcal{H}^{\alpha-1}(m), \\
(1-r^2) f &\longleftrightarrow \mathcal{S} \widehat{F}.
\end{aligned}
\end{equation}
Notice that this operator actually \emph{lowers} $\alpha$.

\begin{figure}[H]
\centering
\includegraphics[width=0.8\linewidth]{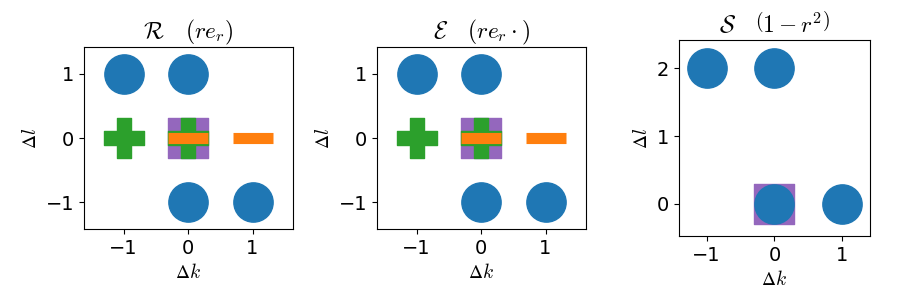}
\caption{Sparse mode coupling for the spherical radial operators.  Unlike the mode coupling for differential operators, we see the radial
    operators increase the $l$ index of a given mode.  This means we must truncate the output of these operators to match the size
    of our truncated expansions.  This in turn implies that, though the differential operators are exactly implemented in our numerical
    scheme, we incur error in the highest modes with radial operations.}
\label{fig:sparsity_radial}
\end{figure}

Figure~\ref{fig:sparsity_radial} demonstrates the mode coupling for each of the radial operators for a single input mode.  
The marker scheme is identical to the differential operator plots in Figure~\ref{fig:sparsity_differential}.

Notice in Figure~\ref{fig:sparsity_radial} that operators may map $k \mapsto k+1$ without lowering the vertical
degree by two.  This is caused by the $r$ multiplication that raises the total polynomial degree; it is unavoidable
that we map out of our basis set here.  Fortunately the consequences of truncation error aren't as strong as when using an incorrect truncation
strategy for the differential operators.

\subsubsection{Conversion}

We must take some care when projecting PDEs onto spherindrical basis modes.  Specifically, the formulation requires term-by-term uniformity in 
${\left( \sigma, \alpha \right)}$ indices which is not automatically guaranteed.
We thus define the family of conversion operators that maps between the Hilbert spaces.  These are identity operators that embed
a field in $\mathcal{H}^{\alpha}(m)$ into $\mathcal{H}^{\alpha + 1}(m)$.  Namely we have
\begin{equation}
\begin{aligned}
\mathcal{I}_{\alpha} &: \mathcal{H}^{\alpha}(m) \to \mathcal{H}^{\alpha+1}(m), \\
f &\longleftrightarrow \mathcal{I}_{\alpha} \widehat{F}.
\end{aligned}
\end{equation}
The family of Hilbert spaces therefore has the nested structure $\mathcal{H}^{\alpha}(m) \subset \mathcal{H}^{\alpha+1}(m)$.

Figure~\ref{fig:conversion-schematic} demonstrates how the embedding operator $\mathcal{I}_{\alpha}$ and multiplication by 
${ \left( 1 - r^{2} \right) }$ act on regularity spaces.  These operators map scalars to scalars and hence $\sigma$ remains
identically zero under their action. 

\begin{figure}[H]
\centering
\begin{tikzpicture}
    \node (A) at (0,0) {$f \in \mathcal{H}^{\alpha}(m)$};
    \node (C) at (4, 0) {$f \in \mathcal{H}^{\alpha+1}(m) $};
    \node (E) at (8,0) {$\left( 1 - r^{2} \right) f \in \mathcal{H}^{\alpha}(m) $};
    \path[->] (A) edge                node[above] {$\mathcal{I}_{\alpha}$} (C);
    \path[->] (C) edge                node[above] {$\mathcal{S}$} (E);
\end{tikzpicture}
\caption{Action of the $\mathcal{I}_{\alpha}$ and $\mathcal{S}$ operators on a scalar field $f$
    with respect to regularity spaces.   Since all fields are scalars the spin weight $\sigma$ is always zero.
    The embedding operator $\mathcal{I}_{\alpha}$ increments the $\alpha$ index, while multiplication by ${ \left(1 - r^{2} \right)}$
    lowers $\alpha$.  $\mathcal{S}$ is the only $\alpha$-lowering operator.}
\label{fig:conversion-schematic}
\end{figure}

\begin{figure}[H]
\centering
\includegraphics[width=0.28\linewidth]{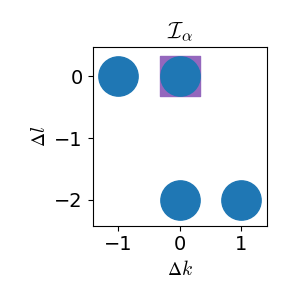}
\caption{Sparse mode coupling for the $\alpha$-conversion operator.  Since the $k$ index is always incremented in tandem with an $l$ decrement of two,
    the conversion operator is exact when triangularly truncating the series expansions.}
\label{fig:sparsity_conversion}
\end{figure}

Figure~\ref{fig:sparsity_conversion} shows the mode coupling for the $\alpha$-conversion operator for a single input mode.  
Here we note that $k$ always increments in tandem with two $l$ decrements and hence the conversion operator is closed in our basis set.

\subsection{Example: Spherical Bessel's Equation}
To utilize the basis for numerical computation we truncate the field expansions up to a specified number of vertical and radial modes $L_{\text{max}}$ and
$N_{\text{max}}$, respectively.  In the eigenvalue problems studied below the systems are decoupled in azimuthal mode so we need only consider a single $m$
at a time.  We take the maximum radial degree to depend on the vertical degree, namely ${ N(l) = N_{\text{max}} - \floor{\frac{l}{2}} }$, in order to maintain
a constant total polynomial degree; this triangular truncation will be evident in the sparsity diagrams below.

After choice of truncation we construct the matrix operators acting on the basis functions.  The operators, defined in Section~\ref{sec:operators} with matrix entries
given explicitly in Appendix A, are linear maps between
basis elements.  We again emphasize that, to achieve maximal sparsity, the operators not only act between $l$ and $k$ indices but also modulate the $\sigma$ and $\alpha$
parameters of the basis.  Note that we are free to select $\alpha$ for each field independently.  Proper choice can improve the
numerical conditioning of the resulting system.

We must take care to ensure each equation is in terms of a single ${ (\sigma, \alpha) }$ index.  This is achieved by converting the
parameters as necessary.  For example, the scalar Laplacian operator maps $\alpha \mapsto \alpha+2$ due to its two derivatives.  To discretize the spherical
Bessel equation
\begin{equation}
\nabla^2 f + \kappa^2 f = 0
\end{equation}
we utilize the $\alpha$ conversion operator ${ \mathcal{I}_{\alpha} : \alpha \to \alpha+1 }$ along with the discretized Laplacian $\mathcal{L}$ and form the matrix equation
\begin{equation}
\mathcal{L} \widehat{F} + \kappa^2 \mathcal{I}_{\alpha}^2 \widehat{F} = 0.
\end{equation}
This defines a generalized eigenproblem for the eigenvalues $\kappa^2$.  What remains is to specify the boundary at ${ r = 1 }$.

\begin{figure}[H]
\centering
\includegraphics[width=.8\linewidth]{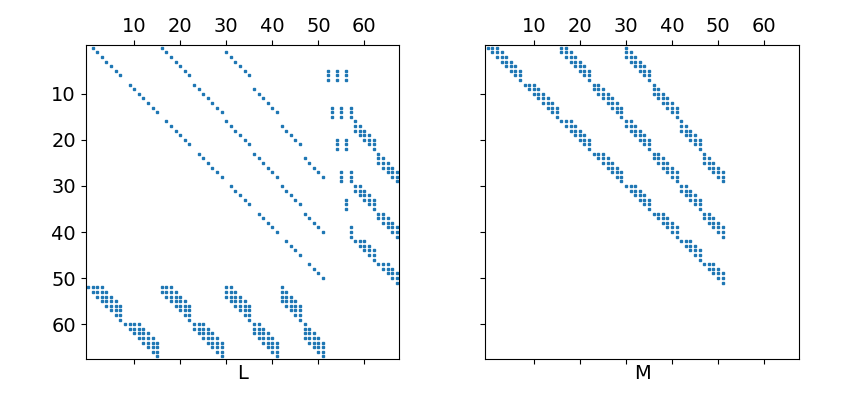}
\caption{Sparsity plot for the Bessel problem.  Tau lines are the dense rows at the bottom of the $L$ matrix.
    Notice they are partitioned into even and odd $l$ indices.  The tau projection coefficients are the sparse
    columns on the far right side of the $L$ matrix}
\label{fig:bessel_spy_tau}
\end{figure}

To implement the boundary condition ${ f(r=1) = 0 }$ we employ two different approaches - the tau method
\cite{lanczos1938trigonometric,gottlieb1977numerical}
and Galerkin recombination \cite{Boyd_2001,shen2011spectral,olver2013fast}.
For the tau method we add extra equations to the system corresponding to evaluation of our field at the boundary.
Because of the coupled spatial directions this produces a number of dense rows that we append to our matrix system.
We then add tau polynomials to the system - extra degrees of freedom that make the system square and soluble.  Choice of tau
polynomial strongly influences the solution; for an excellent summary on the method see Appendix B from \cite{Vasil_Lecoanet_Burns_Oishi_Brown_2019}.
Denoting the boundary evaluation operator $\mathcal{B}$ and the tau polynomial projections to the proper basis $\mathcal{P}$
we form the augmented system
\begin{equation}
\begin{bmatrix} \mathcal{L} & \mathcal{P} \\ \mathcal{B} & 0 \end{bmatrix} \begin{bmatrix} \widehat{F} \\ \tau \end{bmatrix}
    + \kappa^2 \begin{bmatrix} \mathcal{I}_{\alpha}^2 & 0 \\ 0 & 0 \end{bmatrix} \begin{bmatrix} \widehat{F} \\ \tau \end{bmatrix}
    = 0.
\label{eqn:bessel_tau}
\end{equation}
The system (\ref{eqn:bessel_tau}) is a generalized eigenvalue problem of the form ${ L X = \lambda M X }$ for which there exist
several sparse solver packages.
We utilize the \textsc{Umfpack} sparse eigensolver to compute the eigenvalues $\kappa^2$ and corresponding eigenvectors.
Figure~\ref{fig:bessel_spy_tau} displays the sparsity structure for this system for $(L_{\text{max}}, N_{\text{max}}) = (8, 8)$.

\begin{figure}[H]
\centering
\includegraphics[width=.8\linewidth]{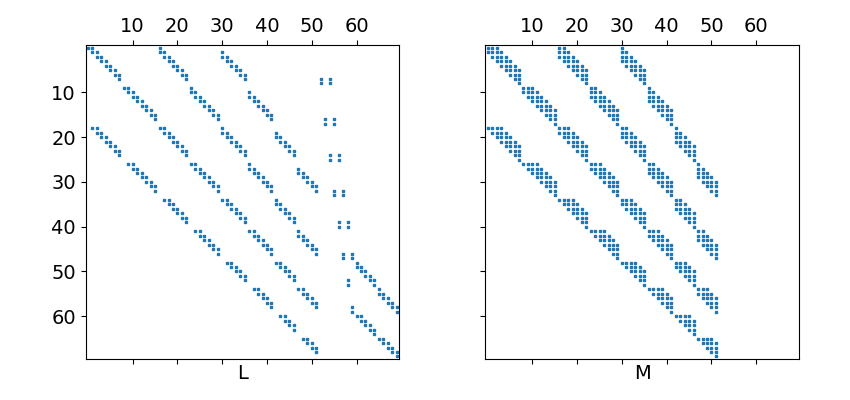}
\caption{Sparsity plot for the Bessel problem using Galerkin recombination.  Note the increase in bandwidth with respect to the tau approach,
    and the absence of dense tau lines across the bottom of the system.  In addition to increasing the bandwidth, the Galerkin recombination
    increases the number of equations in the system.  We are forced to compensate with an equal number of additional tau polynomials, as seen on
    the right hand side of the $L$ and $M$ matrices}
\label{fig:bessel_spy_galerkin}
\end{figure}

The dense boundary rows of the tau method destroy sparsity of the solve and can potentially make the eigensolve quite ill-conditioned,
leading to spurious eigenvalues \cite{gardner1989modified,dawkins1998origin}.
To circumvent this we apply Galerkin recombination of our basis functions to produce a new basis that automatically satisfies the boundary conditions.  
We thus define a family of basis functions that vanish on the boundary: 
\begin{equation}
 \Phi_{m,l,k}^{\sigma, \alpha} \triangleq (1 - r^2) \Psi_{m,l,k}^{\sigma, \alpha+1}.
\end{equation}  
Multiplication by ${ (1-r^2) }$ is a sparse operator $\mathcal{S}$ on basis functions that maps ${ \alpha \mapsto \alpha-1 }$.
We then use the change of variables $f = \mathcal{S} g$ and form the differential equation
\begin{equation}
\nabla^2 \mathcal{S} g + \kappa^2 \mathcal{S} g = 0.
\label{eqn:bessel_galerkin}
\end{equation}
The $\mathcal{S}$ operator increases the maximum polynomial degree: ${ \mathcal{S} \left( L_{\text{max}}, N_{\text{max}} \right) \mapsto \left( L_{\text{max}}+2, N_{\text{max}}+1 \right) }$.
When discretized, the Galerkin system (\ref{eqn:bessel_galerkin}) is therefore underdetermined; one standard remedy is to truncate the output of the $\mathcal{S}$ operator to the input degree.
Equation truncation is equivalent to projecting tau polynomials onto the out-of-range modes.  We elect this more general approach of projecting tau
polynomials to make the system square and soluble, using the $\mathcal{I}_{\alpha}$ operator to project modes from the $\alpha = 1$ basis onto the $\alpha = 2$ equations.
In this case there are no dense boundary rows - just sparse tau projection columns - and so the final system takes the form
\begin{equation}
\begin{bmatrix} \mathcal{L} \mathcal{S} & \mathcal{P} \end{bmatrix} \begin{bmatrix} \widehat{G} \\ \tau \end{bmatrix}
    + \kappa^2 \begin{bmatrix} \mathcal{I}_{\alpha}^2 \mathcal{S} & 0 \end{bmatrix} \begin{bmatrix} \widehat{G} \\ \tau \end{bmatrix}
    = 0.
\end{equation}
Once we have found the coefficients $\widehat{G}$ we then compute the solution $f$ that satisfies the boundary
via ${ \widehat{F} = \mathcal{S} \widehat{G} }$.  Figure~\ref{fig:bessel_spy_galerkin} shows the sparsity diagram for the Galerkin system.

The spherical Bessel’s equation is one with rotational symmetry. Furthermore, spherical harmonics decouple the radial and angular solutions to the spherical Bessel equation.
As such, the spherindrical basis is a poor choice to solve this problem, requiring a two-dimensional solve rather than
the one-dimensional radial problem when using spherical harmonics.  This section instead demonstrates the building blocks required
for a numerical recipe using the spherindrical basis.  More complicated problems lend themselves quite well to the
present work as we demonstrate in the next section.

\subsection{Boundary Evaluation and Tau Polynomials}
To implement a boundary condition at ${ r = 1 }$ using tau polynomials we utilize the boundary evaluation operator $\mathcal{B}$,
which converts the expansion coefficients to an expression for their value on the boundary.  As shown in Appendix A, the operator
decomposes into equations for the even and odd $l$ coefficients.  To set the field equal to zero on the boundary we append these
equations to the discretized matrix system.  We find there are exactly ${ 2 N_{\text{max}} }$
boundary constraint equations in order to set a triangularly truncated field to zero at ${ r = 1 }$.  
This tells us how many tau polynomials are needed to make the system square and soluble once again.

The scalar Laplacian operator in spherindrical coordinates has two coupled spatial derivatives.  After the dust settles these derivatives 
require tau coefficients in the highest radial degree, ${ N_{\text{max}}-1 }$, and highest two vertical degrees, ${ L_{\text{max}}-2 }$ and ${ L_{\text{max}}-1 }$.
We then use our operators to project candidate tau polynomials onto the equation space.  The common choices are the identity operator and various powers of the $\mathcal{I}_{\alpha}$
conversion operator.  From these operators we slice columns corresponding the highest radial and two highest vertical modes then block-append
these to the matrix system.  We stress that tau polynomial choice must be consistent in both spin weight $\sigma$ and numerical index $\alpha$
for the equations where they appear.

Galerkin recombination increases the maximum radial degree by one and the maximum vertical degree by two.  This means we here must also append
tau coefficients to render the system square.  Choice proceeds exactly as described above except with 
${ (L_{\text{max}}, N_{\text{max}}) \mapsto (L_{\text{max}}+2, N_{\text{max}}+1) }$.

\section{Gyroscopic Eigenvalue Test Problems}\label{sec:eigenproblems}
We now solve three eigenvalue problems that demonstrate the efficiency of the basis for the dynamics within a rotating sphere. 
The inertial waves problem requires setting the spherical radial component of the velocity field to zero at ${ r = 1 }$.  This prohibits the use
of Galerkin recombination and demonstrates both the benefits of the basis choice and the numerical ill-conditioning of the
dense boundary operator for large $L_{\text{max}}$.  The damped inertial wave problem adds viscosity to the fluid allowing us
to specify no-slip boundary conditions ${ \vec{u} = 0 }$ at ${ r = 1 }$.  For this problem we demonstrate the effectiveness
of the Galerkin basis.  Our final problem is the linear onset of rotating thermal convection in a ball of fluid.  Despite the Ekman boundary layers
present we show that the basis represents the critical modes in a much sparser expansion than spherical harmonics can attain.

\subsection{Inviscid Inertial Waves}
We solve the inviscid inertial waves problem with velocity and pressure fields $\vec{u}$ and $p$, respectively.  Non-dimensionally this is given by
\begin{equation}
\begin{aligned}
i \lambda \vec{u} + 2 \vec{\hat{e}}_z \vec{\times} \vec{u} &= -\del p \\
\del \cdot \vec{u} &= 0
\end{aligned}
\end{equation}
with impenetrable boundary condition ${ \vec{e}_r \cdot \vec{u} = 0 }$ at ${ r = 1 }$.  The spinor basis (\ref{eqn:spinor_basis}) diagonalizes the Coriolis force
but the boundary condition couples all three velocity components.  Greenspan \cite{Greenspan_1968} provides analytic solutions for both the eigenfrequencies and the associated
pressure modes, thus providing an explicit strategy for testing the accuracy of the spherindrical approach.

To discretize the system we represent pressure with $\alpha = 0$ and velocity with $\alpha = 1$.  These choices bypass any need for conversion
operators in the problem.  In addition we select our state vector 
\begin{equation}
X = \begin{bmatrix} i \widehat{U}^{+} \\ i \widehat{U}^{-} \\ i \widehat{U}^{0} \\ \widehat{P} \end{bmatrix},
\end{equation}
which yields the purely real matrix system, ${ L X = \lambda M X }$.  
Denoting the identity operator $\mathcal{I}$ and boundary evaluation operator $\mathcal{B}$ form the discretized matrix system
\begin{equation}
\lambda \begin{bmatrix} \mathcal{I} & 0 & 0 & 0 & 0 \\ 0 & \mathcal{I} & 0 & 0 & 0 \\ 0 & 0 & \mathcal{I} & 0 & 0 \\ 0 & 0 & 0 & 0 & 0 \\ 0 & 0 & 0 & 0 & 0 \end{bmatrix}
        \begin{bmatrix} i \widehat{U}^{+} \\ i \widehat{U}^{-} \\ i \widehat{U}^{0} \\ \widehat{P} \\ \tau \end{bmatrix}
    = \begin{bmatrix} -2 \mathcal{I} &   0 & 0 & -\mathcal{G}^{+} & \mathcal{P}_{+} \\
                         0 & 2 \mathcal{I} & 0 & -\mathcal{G}^{-} & 0 \\
                         0 & 0 & 0 & -\mathcal{G}^{0} & 0 \\
                         \mathcal{D}^{-} & \mathcal{D}^{+} & \mathcal{D}^{0} & 0 & 0 \\
                         \mathcal{B} \mathcal{E}^{-} & \mathcal{B} \mathcal{E}^{+} & \mathcal{B} \mathcal{E}^{0} & 0 & 0
        \end{bmatrix}
        \begin{bmatrix} i \widehat{U}^{+} \\ i \widehat{U}^{-} \\ i \widehat{U}^{0} \\ \widehat{P} \\ \tau \end{bmatrix}.
\end{equation}
Observe from the matrix system that the only contribution to the $\widehat{U}^{0}$ component of the velocity is $\mathcal{G}^{0} \widehat{P}$.
The operator ${ \mathcal{G}^{0} = \vec{\hat{e}}_{0} \cdot \del }$ always lowers the $l$ index by one.  We improve the conditioning of our
problem by truncating the vertical velocity to maximum vertical degree ${ L_{\text{max}}-1 }$.

\begin{figure}[H]
\centering
\includegraphics[width=.8\linewidth]{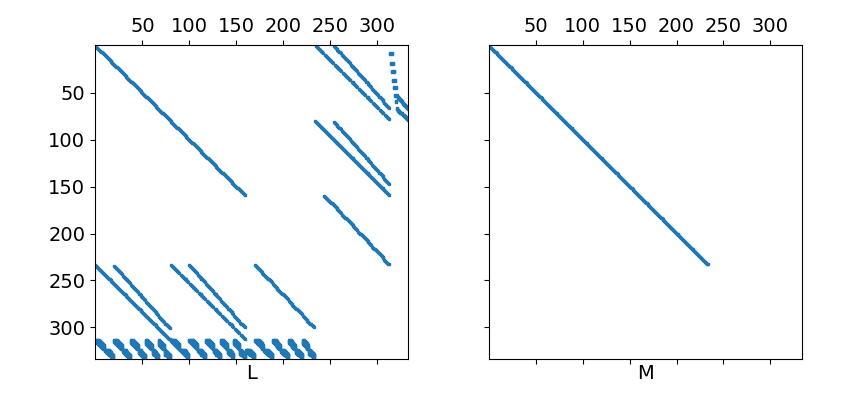}
\caption{Sparsity diagram for the inviscid inertial waves problem for $L_{\text{max}} = N_{\text{max}} = 10$}
\label{fig:inertial_waves_spy}
\end{figure}

Note the problem is first-order in derivatives and hence, by necessity only requires imposition of impenetrable mechanical boundary conditions. 
This means we use a single tau projection operator
\begin{equation}
\mathcal{P}_{+} = \mathcal{I}_{\alpha}[:],
\end{equation}
where the $[:]$ slicing indicates we are taking the projection onto the final radial and final two vertical modes of the expansion
of the $\widehat{U}^{+}$ variable.  The sparsity plot is shown in Figure~\ref{fig:inertial_waves_spy}.

\begin{figure}[H]
\centering
\parbox{\linewidth}{
\includegraphics[width=0.49\linewidth]{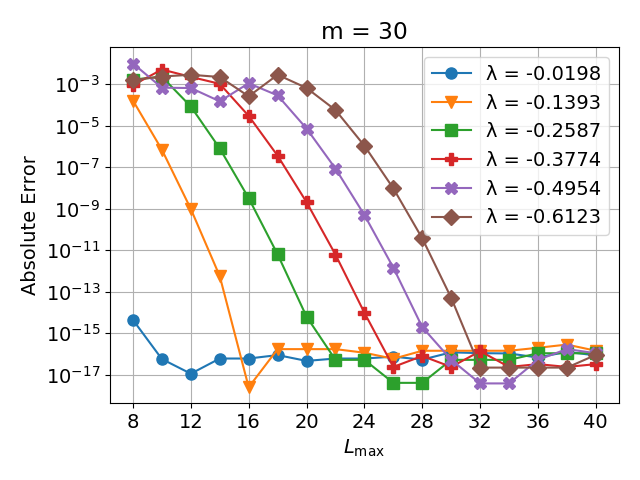}
\includegraphics[width=0.49\linewidth]{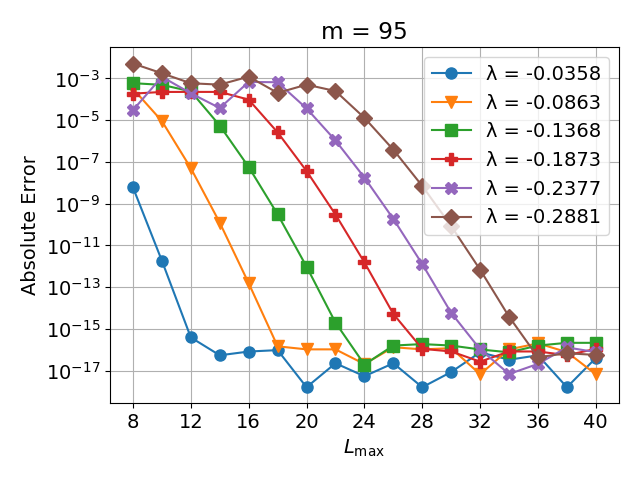}
}
\caption{Inertial wave eigenvalue error as a function of $L_{\text{max}}$ in the $m = 30$ (left) and $m = 95$ (right)
    modes shown in Figures~\ref{fig:inertial_wave_modes_m30} and \ref{fig:inertial_wave_modes_m95}.  The $m = 30$
    plot uses constant radial resolution $N_{\text{max}} = 24$ while the $m = 95$ plot uses $N_{\text{max}} = 32$
    in order to capture all the horizontal oscillations for each mode.  Error is measured with the absolute difference of
    the eigenvalue computed with the eigensolve and Greenspan's analytic results.}
\label{fig:inertial_wave_evalue_error}
\end{figure}

\begin{figure}[H]
\centering
\includegraphics[width=.95\linewidth]{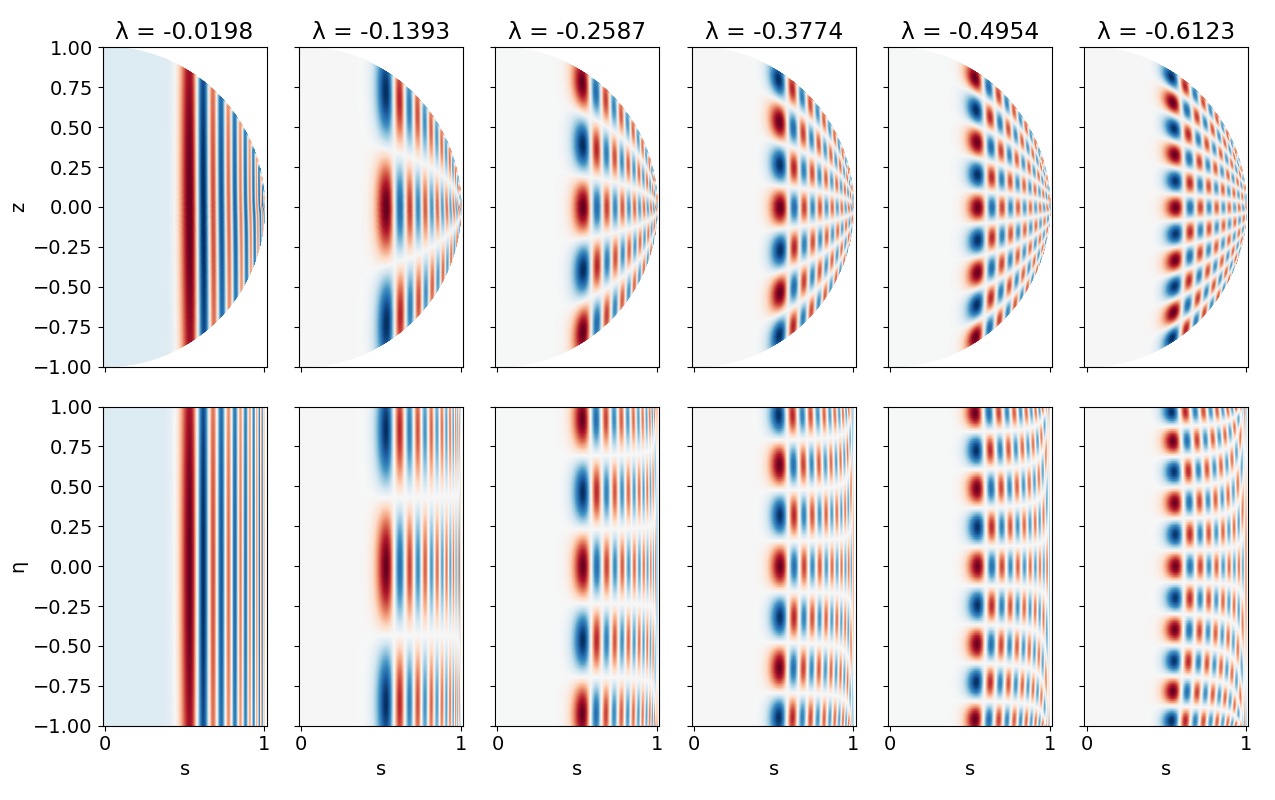}
\caption{Inertial wave eigenfunctions for a single radial degree, $m = 30$.
The upper plots are shown in $(s,z)$ coordinates while the lower plots use $(s, \eta)$ coordinates}
\label{fig:inertial_wave_modes_m30}
\end{figure}

\begin{figure}[H]
\centering
\includegraphics[width=.95\linewidth]{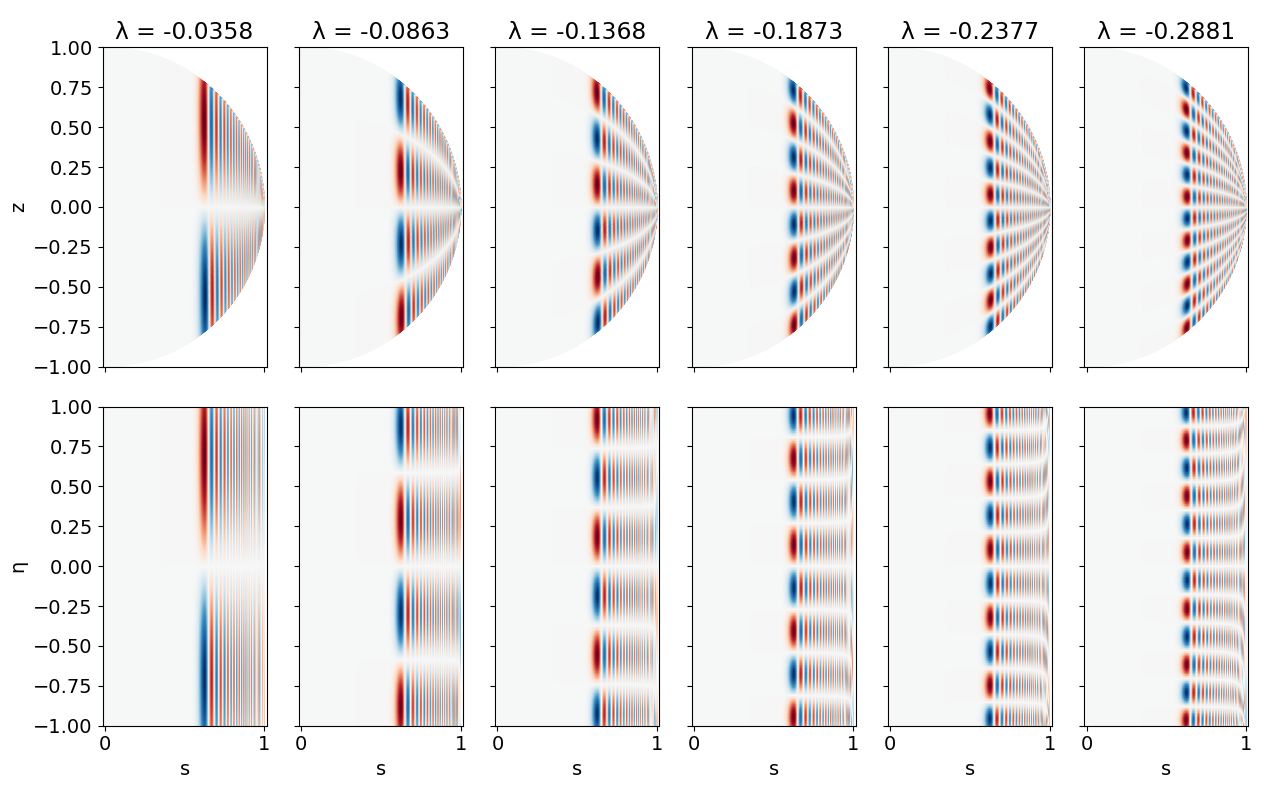}
\caption{Inertial wave eigenfunctions for a single radial degree, $m = 95$.
The upper plots are shown in $(s,z)$ coordinates while the lower plots use $(s, \eta)$ coordinates}
\label{fig:inertial_wave_modes_m95}
\end{figure}

Greenspan writes the solutions in cylindrical coordinates, demonstrating the gyroscopic tendency of rotating fluids to align axially.  
All eigenvalues are real and contained in the interval ${ [-2, 2] }$ and are asymmetric about the imaginary axis.
Figure~\ref{fig:inertial_wave_evalue_error} shows the absolute error between a selection of analytic eigenvalues and their numerical solution for ${m = 30}$ and ${m = 95}$.
We achieve machine precision for the modes with low degree vertical oscillations with very modest resolution requirements.
Figures~\ref{fig:inertial_wave_modes_m30} and \ref{fig:inertial_wave_modes_m95} plot meridional slices of the corresponding modes.  
Flow structures of modes with growth rates near zero align vertically, with small variation in the vertical direction compared to the horizontal.

Resolving these structures with spherical harmonics requires a larger expansion in harmonic degree than needed for the spherinder basis.
For the $m = 95$ solutions, the spherinder basis with ${ \left( L_{\text{max}}, N_{\text{max}} \right) = \left( 36, 32 \right) }$ resolves all target modes to machine precision, using 3,433 total degrees of freedom.
Spherical harmonics require angular degree $L_{\text{max}} = 155$ and spherical radial degree $N_{\text{max}} = 77$ to resolve these modes to the same precision, but instead require
4,269 degrees of freedom.  These performance gains become more pronounced with increased gyroscopic alignment.

\subsection{Damped Inertial Waves}
The damped inertial wave equations model the exponential decay in time of inertial modes for fluids with viscosity.  We solve
\begin{equation}
\begin{aligned}
\partial_t \vec{u} + 2 \vec{\hat{e}}_z \vec{\times} \vec{u} &= -\del p + \textrm{E } \del^2 \vec{u} \\
\del \cdot \vec{u} &= 0
\end{aligned}
\end{equation}
subject to no-slip boundary conditions ${ \vec{u} = 0 }$ at ${ r = 1 }$.  Here we use the Ekman number
definition ${ E = \frac{ \nu }{ \Omega R^2 } }$, where $\nu$ is the kinematic viscosity, $\Omega$ is 
the rotation rate and $R = 1$ is the radius of the sphere.  The Ekman number denotes the relative importance of viscous diffusion to the fluid’s inertial acceleration force. 
We make the ansatz ${ \partial_t \mapsto \lambda }$ (i.e. exponential time dependence) and
solve for the eigenvalues $\lambda$ of the discretized system.  To implement the boundary condition we utilize Galerkin recombination for all three
components of the velocity.

The matrix system takes the form
\begin{equation}
L X = \lambda M X
\end{equation}
where
\begin{equation}
X =  \begin{bmatrix} \widehat{U}^{+} \\ \widehat{U}^{-} \\ \widehat{U}^{0} \\ \widehat{P} \\ \tau \end{bmatrix}, \hspace{8ex}
M = \begin{bmatrix} 
    \mathcal{I}_{\alpha}^2 \mathcal{S}_{+} & 0 & 0 & 0 & 0 \\ 
    0 & \mathcal{I}_{\alpha}^2 \mathcal{S}_{-} & 0 & 0 & 0 \\ 
    0 & 0 & \mathcal{I}_{\alpha}^2 \mathcal{S}_{0} & 0 & 0 \\ 
    0 & 0 & 0 & 0 & 0
\end{bmatrix}
\end{equation}
and
\begin{equation}
L = \begin{bmatrix} 
    \left( \textrm{E} \mathcal{L}_{+} - 2 i \mathcal{I}_{\alpha}^2 \right) \mathcal{S}_{+} & 0 & 0 & -\mathcal{G}^{+} & \mathcal{P}_{+} \\ 
     0 & \left( \textrm{E} \mathcal{L}_{-} + 2 i \mathcal{I}_{\alpha}^2 \right) \mathcal{S}_{-} & 0 & -\mathcal{G}^{-} & \mathcal{P}_{-} \\ 
     0 & 0 & \textrm{E} \mathcal{L}_{0} \mathcal{S}_{0} & -\mathcal{G}^{0} & \mathcal{P}_{0} \\ 
     \mathcal{D}^{-} \mathcal{S}_{+} & \mathcal{D}^{+} \mathcal{S}_{-} & \mathcal{D}^{0} \mathcal{S}_{0} & 0 & \mathcal{P}_{\text{div}}.
\end{bmatrix}
\end{equation}
The momentum equations live in $\mathcal{H}(\sigma,3)$ while the divergence equations lives in $\mathcal{H}(0,2)$.
The cascaded conversion operators in the $M$ matrix and in the Coriolis terms make the momentum equations consistent with $\alpha = 3$.
We choose $\alpha = 1$ for the velocity and $\alpha = 2$ for the pressure field so that we can avoid converting the pressure
gradients to $\alpha = 3$.  Shifting all fields by a fixed $\alpha$ offset retains the matrix structure but can improve conditioning.
As in the inertial waves problem we truncate vertical velocity to a maximum vertical degree ${ L_{\text{max}}-1 }$ since only the $z$
derivative of the pressure contributes.

\begin{figure}[H]
\centering
\includegraphics[width=.49\linewidth]{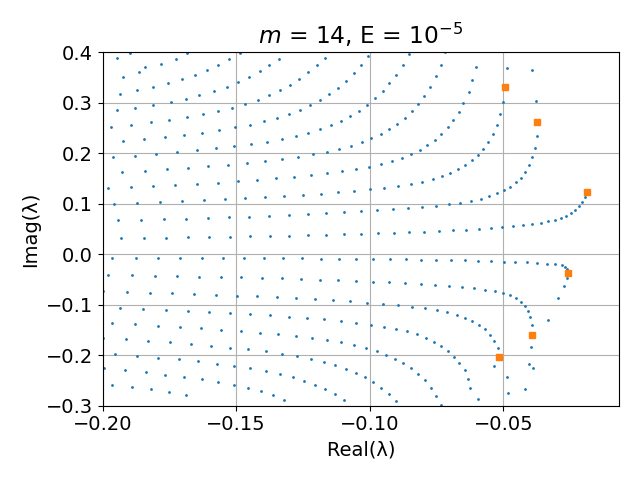}
\includegraphics[width=.49\linewidth]{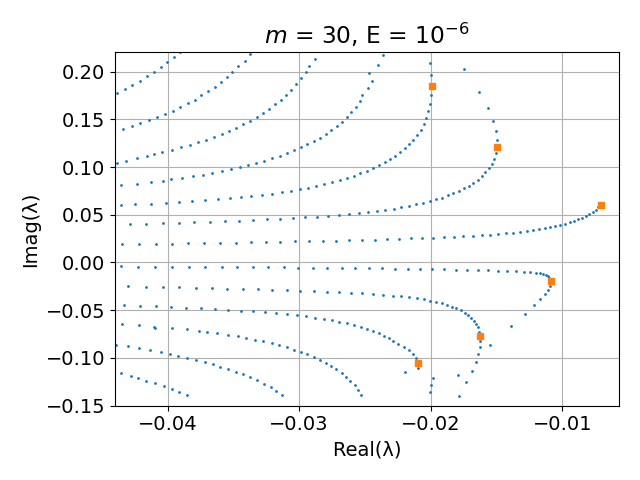}
\caption{Damped inertial wave eigenvalues, $m = 14$, $\textrm{E} = 10^{-5}$ (left) and $m = 30$, $\textrm{E} = 10^{-6}$ (right).  
    Each eigenvalue curve represents a different vertical degree.  Following the curve from right to left corresponds with an increase in
    radial degree.  The orange square symbols denote modes selected for plotting below.}
\label{fig:hydrodynamics_evalues}
\end{figure}

\begin{figure}[H]
\centering
\includegraphics[width=\linewidth]{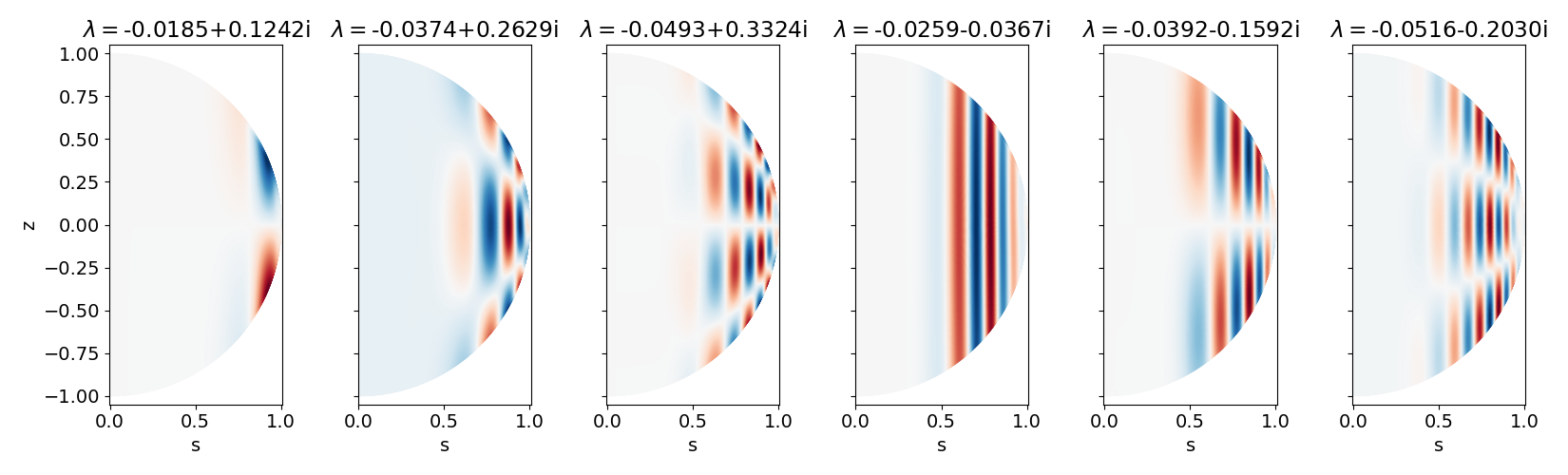}
\caption{Damped inertial pressure modes, $m = 14$, $\textrm{E} = 10^{-5}$.  Displayed is the least damped mode of the six most critical eigenvalue curves
    indicated with orange squares in Figure~\ref{fig:hydrodynamics_evalues} (left)}
\label{fig:hydrodynamics_evectors_m14}
\end{figure}

\begin{figure}[H]
\centering
\includegraphics[width=\linewidth]{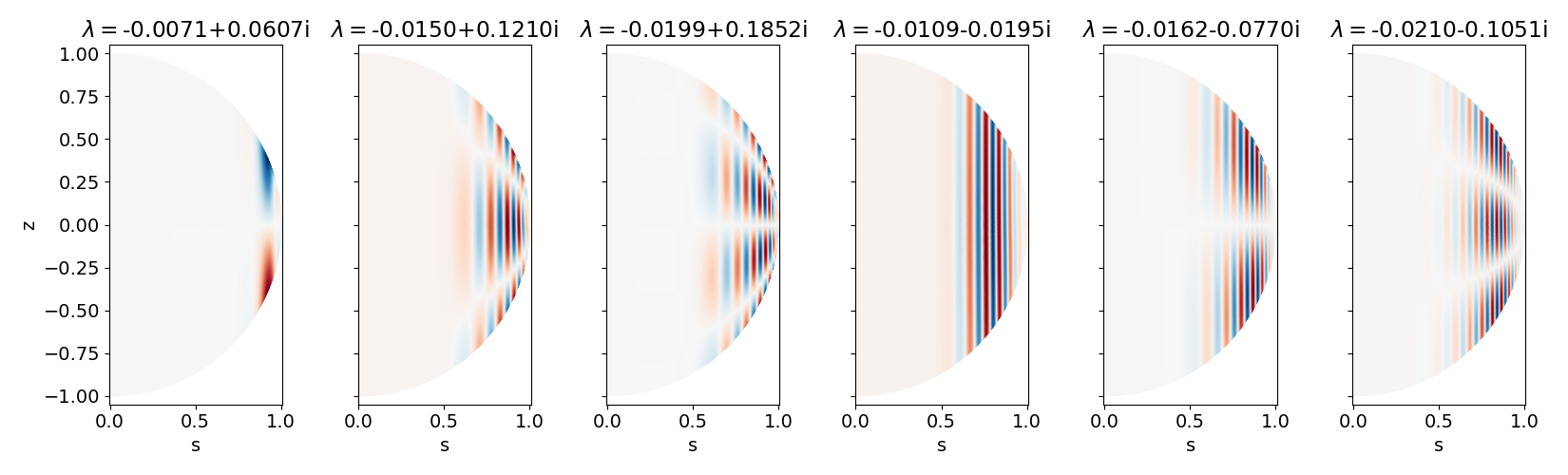}
\caption{Damped inertial pressure modes, $m = 30$, $\textrm{E} = 10^{-6}$.  Displayed is the least damped mode of the six most critical eigenvalue curves
    indicated with orange squares in Figure~\ref{fig:hydrodynamics_evalues} (right)}
\label{fig:hydrodynamics_evectors_m30}
\end{figure}

Figure~\ref{fig:hydrodynamics_evalues} displays the computed eigenvalues for $m = 14$, $\textrm{E} = 10^{-5}$ and $m = 30$, $\textrm{E} = 10^{-6}$.
The orange squares mark the eigenfrequencies for which we plot the pressure mode in Figure~\ref{fig:hydrodynamics_evectors_m14} (${ m = 14 }$, ${ \textrm{E} = 10^{-5} }$)
and Figure~\ref{fig:hydrodynamics_evectors_m30} (${ m = 30 }$, ${ \textrm{E} = 10^{-6} }$).  These minimally damped modes display rapid variation
in the cylindrical radial direction but slow variation in the vertical.  The spherindrical basis excels at representing these vertically aligned flows.

\subsubsection{Comparison with Spherical Harmonic Representation}
\begin{figure}[H]
\centering
\includegraphics[width=0.49\linewidth]{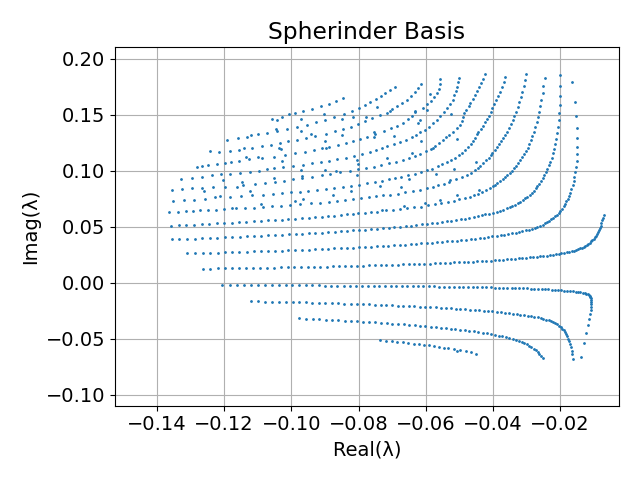}
\includegraphics[width=0.49\linewidth]{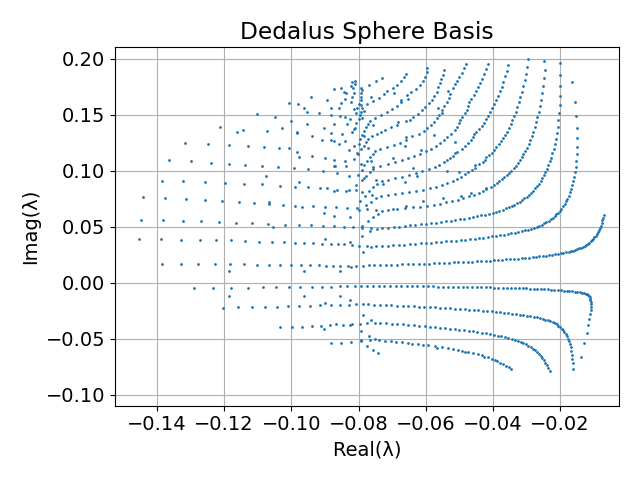}
\caption{Damped inertial wave eigenvalues for ${ m = 30 }$, ${ \textrm{E} = 10^{-6} }$.  Shown are the thousand modes nearest the least damped one for the spherinder basis (left) and \textsc{Dedalus} sphere basis (right).}
\label{fig:hydrodynamics_evalues-comparison}
\end{figure}

We now compare the spherinder basis resolution capabilities to spherical harmonics code.  We implement the damped inertial
waves problem using the \textsc{Dedalus} software package \cite{Burns_Vasil_Oishi_Lecoanet_Brown_2020, Vasil_Lecoanet_Burns_Oishi_Brown_2019, Lecoanet_Vasil_Burns_Brown_Oishi_2019}.
This code uses spin-weighted spherical harmonics and Zernike polynomials in the spherical radial direction to represent tensor fields while explicitly handling coordinate singularities.
To get a feel for the resolution capabilities of the spherinder basis versus that of \textsc{Dedalus} we discretize the
${ m = 30 }$, ${ \textrm{E} = 10^{-6} }$ eigenproblem with comparable number of degrees of freedom between the two bases.  
We use just over 72,000 degrees of freedom with resolution ${ (L_{\text{max}}, N_{\text{max}}) = (80, 240) }$
in the spherinder basis (72,287 total degrees of freedom) and ${ (L_{\text{max}}, N_{\text{max}}) = (280, 146) }$ in the \textsc{Dedalus} sphere basis (72,289 total degrees of freedom).
We then solve for the thousand eigenvalues nearest the most critical one at ${ \lambda \approx -0.0071 + 0.0607i }$.
Figure~\ref{fig:hydrodynamics_evalues-comparison} displays results for the two bases.

Notice the spherinder basis continues the smooth eigenvalue curves farther into the negative reals than the \textsc{Dedalus}
sphere code.  These eigenvalues are extremely oscillatory in the cylindrical radial direction.  Each eigenvalue branch corresponds
to a different vertical degree, and the radial degree increases as the tracks are followed leftwards.  The vertical degree increases
as we depart from the real line.
Both bases struggle in the upper left portion of the diagram.  This region is spatially extremely oscillatory and requires more resolution
to capture the eigenvalues accurately.  Increasing resolution for either basis improves convergence in this region.  Doing so pushes out
into the left-half plane the adequately resolved modes, but the trend remains - the spherinder basis resolves modes nearly 
twice as far into the left-half plane as the \textsc{Dedalus} sphere basis with comparable degrees of freedom.

\subsection{Rotating Thermal Convection}
The final eigenproblem we solve is the linear onset to rotating thermal convection.  In the rapidly rotating limit the critical mode
structure is cylindrical, with rapid variation in both the $s$ and $\phi$ directions but slow variation in the vertical $z$
direction.  The Taylor-Proudman theorem provides a leading order dominant (geostrophic)
balance between the Coriolis and pressure gradient forces which enforces this anisotropy.  The nondimensional perturbation equations take the form
\begin{equation}
\begin{aligned}
\textrm{E} \left( \partial_t - \del^2 \right) \vec{u} + \vec{\hat{e}}_z \vec{\times} \vec{u} &= -\del p + \textrm{Ra } \vartheta \vec{r} \\
\del \cdot \vec{u} &= 0 \\
\left( \textrm{Pr } \partial_t - \nabla^2 \right) \vartheta &= \textrm{Pr } \vec{u} \cdot \vec{r},
\end{aligned}
\end{equation}
where $\vartheta$ is the deviation from the base temperature profile ${ \frac{1}{2} (1-r^2) }$.  For details of the nondimensionalization see \cite{Marti_Schaeffer_Hollerbach_2014}.

We impose the no-slip boundary condition ${ \vec{u} = 0 }$
and the fixed temperature condition ${ \vartheta = 0 }$ at ${ r = 1 }$, which suggests Galerkin recombination to specify the boundary data.

\begin{figure}[H]
\centering
\includegraphics[width=\linewidth]{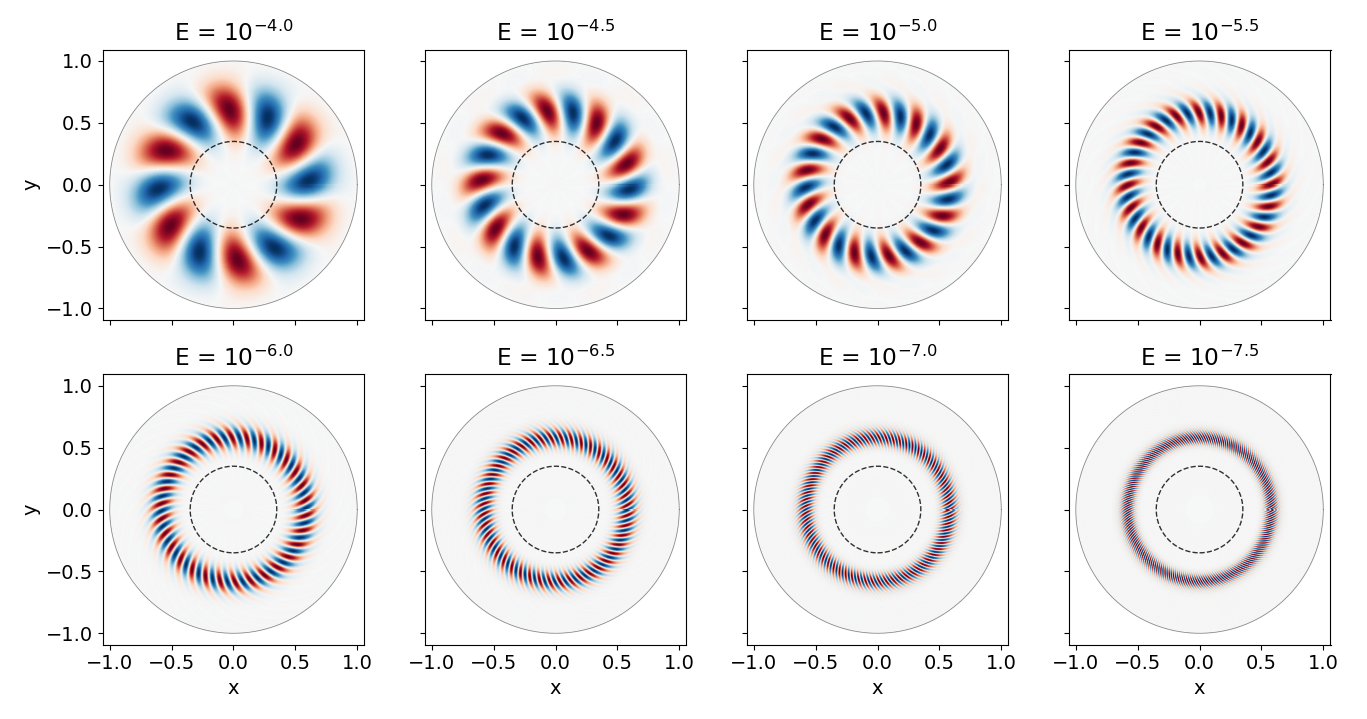}
\caption{Equatorial slices of the critical temperature field for the rotating thermal convection problem computed with the spherinder basis, plotted for various Ekman numbers.  We superimpose a fictitious inner core of
    radius $r_{i}/r_{o} = 0.35$ to provide visual confirmation for convergence of critical Rayleigh number to the shell results $\widetilde{\textrm{Ra}}_{M}$.  The mode spiraling and thin cylindrical structures as 
    ${E \to 0}$ are apparent in these slices.}
\label{fig:linear_onset-critical_modes-equatorial}
\end{figure}

\begin{figure}[H]
\centering
\includegraphics[width=\linewidth]{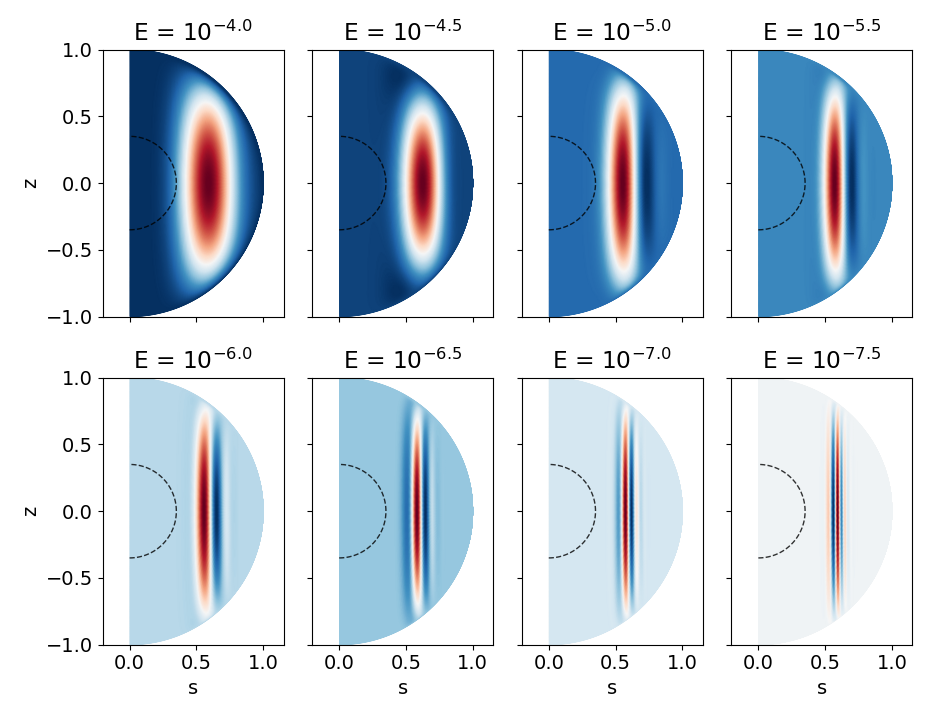}
\caption{Meridional slices of the critical temperature field for the rotating thermal convection problem computed with the spherinder basis, plotted for various Ekman numbers.  We superimpose a fictitious inner core of
    radius $r_{i}/r_{o} = 0.35$ to provide visual confirmation for convergence of critical Rayleigh number to the shell results $\widetilde{\textrm{Ra}}_{M}$}
\label{fig:linear_onset-critical_modes}
\end{figure}

The critical Rayleigh number $\textrm{Ra}_{c}$ is that which sends the real part of the least damped eigenvalue to zero.  At this thermal forcing
the least damped eigenvalue $\lambda_{c}$ then takes the form ${ \lambda_{c} = i \omega_{c} }$, where $\omega_{c}$ is the critical frequency.
We compute the critical Rayleigh number for a range of Ekman numbers setting ${ \textrm{Pr} = 1 }$.
Following the scaling of Marti et al \cite{Marti_Calkins_Julien_2016} we take the following definitions for reduced Rayleigh number
and reduced frequency:
\begin{equation}
\begin{aligned}
\widetilde{\textrm{Ra}}_{c} &= \textrm{Ra}_{c} \textrm{E}^{4/3} \hspace{4ex} &&\text{(reduced critical Rayleigh number)} \\
\widetilde{\omega}_{c} &= \omega_{c} \textrm{E}^{2/3} \hspace{4ex} &&\text{(reduced critical frequency)}.
\end{aligned}
\end{equation}
The inner core at radius $r_{i}/r_{o} = 0.35$ in the Marti simulations has little effect on the critical modes, as observed in \cite{Dormy_Soward_Jones_Jault_Cardin_2004},
especially for the rapid rotation regime.  For this reason we provide for comparison $\widetilde{\textrm{Ra}}_{M}$, the critical Rayleigh numbers computed by
Marti et al in Table~\ref{tab:linear_onset}.  The amplitude envelope of the critical mode becomes sharper as we increase rotation rate.  When this envelope is
sharp enough the critical mode dynamics are well-separated from the inner core and hence its influence vanishes.  We show the full sphere critical modes along
with a superimposed inner core in Figures~\ref{fig:linear_onset-critical_modes-equatorial},~\ref{fig:linear_onset-critical_modes} to demonstrate this effect.

\begin{table}[H]
\begin{tabular}{|l|l|l|l|l|}
\hline
\textrm{E} & $m$ & $\widetilde{\omega}_{c}$ & $\widetilde{\textrm{Ra}}_{c}$ & $\widetilde{\textrm{Ra}}_{M}$ \\
\hline
$10^{-4}  $ & 6   & -0.27009 & 5.0151 & 5.1324 \\
$10^{-4.5}$ & 9   & -0.31035 & 4.6581 & 4.6814 \\
$10^{-5}  $ & 13  & -0.33901 & 4.6581 & 4.4665 \\
$10^{-5.5}$ & 20  & -0.36778 & 4.3488 & 4.3484 \\
$10^{-6}  $ & 30  & -0.38930 & 4.2736 & 4.2736 \\
$10^{-6.5}$ & 44  & -0.40439 & 4.2236 & 4.2235 \\
$10^{-7}  $ & 64  & -0.41737 & 4.1903 & 4.1902 \\
$10^{-7.5}$ & 95  & -0.42658 & 4.1677 & 4.1677 \\
\hline
\end{tabular}
\caption{Critical Rayleigh numbers and eigenfrequencies for various Ekman numbers.  Observe the results converge to
    those of the Marti et al despite the presence of an inner core of radius $r_{i}/r_{o} = 0.35$ in their simulations.}
\label{tab:linear_onset}
\end{table}

\section{Conclusions}\label{sec:conclusion}
We introduced the spherindrical coordinate system and corresponding orthogonal basis for numerical computations in the sphere.
The coordinate system breaks rotational spherical symmetry by design to adhere to
gyroscopically aligned dynamics such as those occurring in rotationally constrained geophysical objects. 
This enables a sparser representation of flows impacted by the Taylor-Proudman constraint than can be achieved using
the traditional methodology of spherical harmonics.  This is borne out in the three rotating test problems where,
as the problems generate increasingly small gyroscopic scales, the spherindrical method becomes superior to a
spherical harmonics approach.

A major theme of this work is designing numerical methods for the geometry of the problem.  
Gyroscopic alignment of rapidly rotating fluids motivates our choice to eschew spherical coordinates with their associated orthogonal
polynomials and instead utilize spherindrical coordinates.  Once we specified our coordinate system we sought a class of bases
that conforms to the coordinate singularities.  Proper choice of basis - namely in selection of Jacobi polynomial parameters - yields
a numerical algorithm explicitly free of singularities.  Calculus operators map between bases in the hierarchy in a sparse way.
In this manner we achieve a sparse linear algebraic representation geared for rotating fluids problems that can be solved with a small number of
degrees of freedom.

We demonstrated the feasibility of the gyroscopic approach and opened the door to future detailed investigations of linear stability
analysis and fully nonlinear time-domain simulations.
Generalized eigenvalue problems $L X = \lambda M X$ like those explored in this paper translate directly to time-dependent simulation.
One need only replace $\lambda \mapsto \partial_{t}$ and use a numerical time-stepping algorithm to evolve the system.
For efficiency, nonlinear operators demand grid-space evaluation; spectral transforms are therefore the last required simulation components.
The orthogonal structure of the spherinder basis enables highly accurate transforms through Gauss quadrature.  We will detail these
algorithms in a future publication. 

\bigskip
\noindent \textbf{Acknowledgments: } A.E and K.J. acknowledge support from NSF Grant DMS-2009319.

\begin{appendices}

\section{Vector Calculus for the Spherindrical Basis}
Vasil et al \cite{Vasil_Lecoanet_Burns_Oishi_Brown_2019} define all required Jacobi polynomial operations. 
These fundamental operator definitions generate the matrix coefficients for the calculus operators on our basis
functions.  The Jacobi embedding operators are:
\begin{equation}
\begin{aligned}
\mathcal{I}_{a} P_{n}^{(a,b)}(z) & &&= \hphantom{-}\frac{n+a+b+1}{2n+a+b+1} P_{n}^{(a+1,b)}(z) - \frac{n+b}{2n+a+b+1} P_{n-1}^{(a+1,b)}(z) \\
\mathcal{I}_{a}^{\dagger} P_{n}^{(a,b)}(z) &= (1-z) P_{n}^{(a,b)}(z) &&= -\frac{2(n+1)}{2n+a+b+1} P_{n+1}^{(a-1,b)}(z) + \frac{2(n+a)}{2n+a+b+1} P_{n}^{(a-1,b)}(z) \\
\mathcal{I}_{b} P_{n}^{(a,b)}(z) & &&= \hphantom{-}\frac{n+a+b+1}{2n+a+b+1} P_{n+1}^{(a,b+1)}(z) + \frac{n+a}{2n+a+b+1} P_{n}^{(a,b+1)}(z) \\
\mathcal{I}_{b}^{\dagger} P_{n}^{(a,b)}(z) &= (1+z) P_{n}^{(a,b)}(z) &&= \hphantom{-}\frac{2(n+1)}{2n+a+b+1} P_{n+1}^{(a,b+1)}(z) + \frac{2(n+b)}{2n+a+b+1} P_{n}^{(a,b+1)}(z).
\end{aligned}
\end{equation}
For future use define the coefficients $\gamma_{l}^{(\alpha)}$ and $\delta_{l}^{(\alpha)}$ by the relation
\begin{equation}
P_{l}^{(\alpha,\alpha)}(\eta) = \gamma_{l}^{(\alpha)} P_{l}^{(\alpha+1,\alpha+1)}(\eta) - \delta_{l}^{(\alpha)} P_{l-2}^{(\alpha+1,\alpha+1)}(\eta),
\end{equation}
which are the coefficients of the $(a,b)$ raising operator, $\mathcal{I}_{a} \mathcal{I}_{b}$.
In addition we define $c_{l}^{(\alpha)}$ and $d_{l}^{(\alpha)}$ such that
\begin{equation}
\eta P_{l}^{(\alpha,\alpha)}(\eta) = c_{l}^{(\alpha)} P_{l+1}^{(\alpha,\alpha)}(\eta) + d_{l}^{(\alpha)} P_{l-1}^{(\alpha, \alpha)}(\eta).
\end{equation}
We compute $c_{l}^{(\alpha)}$ and $d_{l}^{(\alpha)}$ from the Jacobi operator
$\mathcal{Z} = \frac{1}{2} \left( \mathcal{I}_{b}^{\dagger} \mathcal{I}_{b} - \mathcal{I}_{a}^{\dagger} \mathcal{I}_{a} \right)$.

The Jacobi differential operators are:
\begin{equation}
\begin{aligned}
\mathcal{D}_{m} P_{n}^{(a,b)}(z)           &= \frac{d}{dz} P_{n}^{(a,b)}(z) &&= \frac{n+a+b+1}{2} P_{n-1}^{(a+1,b+1)}(z) \\
\mathcal{D}_{m}^{\dagger} P_{n}^{(a,b)}(z) &= \left[ (1+z)a - (1-z)b - (1-z^2) \frac{d}{dz} \right] P_{n}^{(a,b)}(z) &&= 2(n+1) P_{n+1}^{(a-1,b-1)}(z) \\
\mathcal{D}_{s} P_{n}^{(a,b)}(z)           &= \left[ b + (1+z) \frac{d}{dz} \right] P_{n}^{(a,b)}(z) &&= (n+b) P_{n}^{(a+1,b-1)}(z) \\
\mathcal{D}_{s}^{\dagger} P_{n}^{(a,b)}(z) &= \left[ a - (1-z) \frac{d}{dz} \right] P_{n}^{(a,b)}(z) &&= (n+a) P_{n}^{(a-1,b+1)}(z) .
\end{aligned}
\end{equation}
Again to aid notation below we define
\begin{equation}
\frac{d}{d\eta} P_{l}^{(\alpha,\alpha)}(\eta) = \beta_{l}^{(\alpha)} P_{l-1}^{(\alpha+1,\alpha+1)},
\end{equation}
which is the super-diagonal term in the $\mathcal{D}_{m}$ operator.

In what follows we expand calculus operators into their various $l$ output components, leaving radial dependence in terms of Jacobi operators.
Results for action on a basis element $\Psi_{m,l,k}^{\sigma,\alpha}$ can always be written as a linear combination
of a few nearby (in $l$ and $k$) basis vectors.  The size of the linear combination determines the sparsity of the
matrix system.  Typical operations map a single mode to between two and four modes in the output space.

\subsection{Scalar Gradient}
\begin{subequations}
\begin{flalign}
&\mathcal{G}^{-}    :\hspace{4ex} \vec{\hat{e}}_{-}^*       \cdot \nabla \Psi_{m,l,k}^{0,\alpha} = \Psi_{m,l,\bigcdot}^{-,\alpha+1} \cdot 2 \gamma_{l}^{(\alpha)} \mathcal{D}_{s} + \Psi_{m,l-2,\bigcdot}^{-,\alpha+1} \cdot 2 \delta_{l}^{(\alpha)} \mathcal{D}_{m}^{\dagger} & \\
&\mathcal{G}^{+}    :\hspace{4ex} \vec{\hat{e}}_{+}^*       \cdot \nabla \Psi_{m,l,k}^{0,\alpha} = \Psi_{m,l,\bigcdot}^{+,\alpha+1} \cdot 2 \gamma_{l}^{(\alpha)} \mathcal{D}_{m} + \Psi_{m,l-2,\bigcdot}^{+,\alpha+1} \cdot 2 \delta_{l}^{(\alpha)} \mathcal{D}_{s}^{\dagger} & \\
&\mathcal{G}^{\zero}:\hspace{4ex} \vec{\hat{e}}_{\zero}^{*} \cdot \nabla \Psi_{m,l,k}^{0,\alpha} = \Psi_{m,l-1,k}^{0,\alpha+1} \cdot \sqrt{2} \beta_{l}^{(\alpha)} &
\end{flalign}
\end{subequations}

\subsection{Vector Divergence}
\begin{subequations}
\begin{flalign}
&\mathcal{D}^{+}    :\hspace{4ex} \nabla \cdot \left(\vec{\hat{e}}_{-}     \Psi_{m,l,k}^{-,\alpha} \right) = \Psi_{m,l,\bigcdot}^{0,\alpha+1} \cdot 2 \gamma_{l}^{(\alpha)} \mathcal{D}_{m} + \Psi_{m,l-2,\bigcdot}^{0,\alpha+1} \cdot 2 \delta_{l}^{(\alpha)} \mathcal{D}_{s}^{\dagger} & \\
&\mathcal{D}^{-}    :\hspace{4ex} \nabla \cdot \left(\vec{\hat{e}}_{+}     \Psi_{m,l,k}^{+,\alpha} \right) = \Psi_{m,l,\bigcdot}^{0,\alpha+1} \cdot 2 \gamma_{l}^{(\alpha)} \mathcal{D}_{s} + \Psi_{m,l-2,\bigcdot}^{0,\alpha+1} \cdot 2 \delta_{l}^{(\alpha)} \mathcal{D}_{m}^{\dagger} & \\
&\mathcal{D}^{\zero}:\hspace{4ex} \nabla \cdot \left(\vec{\hat{e}}_{\zero} \Psi_{m,l,k}^{0,\alpha} \right) = \Psi_{m,l-1,k}^{0,\alpha+1} \cdot \sqrt{2} \beta_{l}^{(\alpha)} &
\end{flalign}
\end{subequations}

\subsection{Curl}
\begin{subequations}
\begin{flalign}
&\mathcal{C}_{-}^{\zero}:\hspace{4ex} \vec{\hat{e}}_{-}^{*}     \cdot \del \times \left( \vec{\hat{e}}_{-} \Psi_{m,l,k}^{-,\alpha} \right) = \Psi_{m,l-1,k}^{-,\alpha+1} \cdot \left( -i \sqrt{2} \beta_{l}^{(\alpha)} \right) & \\
&\mathcal{C}_{\zero}^{+}:\hspace{4ex} \vec{\hat{e}}_{\zero}^{*} \cdot \del \times \left( \vec{\hat{e}}_{-} \Psi_{m,l,k}^{-,\alpha} \right) = \Psi_{m,l,\bigcdot}^{0,\alpha+1} \cdot \left( +2 i \gamma_{l}^{(\alpha)} \mathcal{D}_{m} \right) + \Psi_{m,l-2,\bigcdot}^{0,\alpha+1} \cdot \left( +2 i \delta_{l}^{(\alpha)} \mathcal{D}_{s}^{\dagger} \right) & \\
&\mathcal{C}_{+}^{\zero}:\hspace{4ex} \vec{\hat{e}}_{+}^{*}     \cdot \del \times \left( \vec{\hat{e}}_{+} \Psi_{m,l,k}^{+,\alpha} \right) = \Psi_{m,l-1,k}^{+,\alpha+1} \cdot \left( +i \sqrt{2} \beta_{l}^{(\alpha)} \right) & \\
&\mathcal{C}_{\zero}^{-}:\hspace{4ex} \vec{\hat{e}}_{\zero}^{*} \cdot \del \times \left( \vec{\hat{e}}_{+} \Psi_{m,l,k}^{+,\alpha} \right) = \Psi_{m,l,\bigcdot}^{0,\alpha+1} \cdot \left( -2 i \gamma_{l}^{(\alpha)} \mathcal{D}_{s} \right) + \Psi_{m,l-2,\bigcdot}^{0,\alpha+1} \cdot \left( -2 i \delta_{l}^{(\alpha)} \mathcal{D}_{m}^{\dagger} \right) & \\
&\mathcal{C}_{-}^{-}    :\hspace{4ex} \vec{\hat{e}}_{-}^{*}     \cdot \del \times \left( \vec{\hat{e}}_{\zero} \Psi_{m,l,k}^{0,\alpha} \right) = \Psi_{m,l,\bigcdot}^{-,\alpha+1} \cdot \left( +2 i \gamma_{l}^{(\alpha)} \mathcal{D}_{s} \right) + \Psi_{m,l-2,\bigcdot}^{-,\alpha+1} \cdot \left( +2 i \delta_{l}^{(\alpha)} \mathcal{D}_{m}^{\dagger} \right) & \\
&\mathcal{C}_{+}^{+}    :\hspace{4ex} \vec{\hat{e}}_{+}^{*}     \cdot \del \times \left( \vec{\hat{e}}_{\zero} \Psi_{m,l,k}^{0,\alpha} \right) = \Psi_{m,l,\bigcdot}^{+,\alpha+1} \cdot \left( -2 i \gamma_{l}^{(\alpha)} \mathcal{D}_{m} \right) + \Psi_{m,l-2,\bigcdot}^{+,\alpha+1} \cdot \left( -2 i \delta_{l}^{(\alpha)} \mathcal{D}_{s}^{\dagger} \right) &
\end{flalign}
\end{subequations}

\subsection{Multiplication by $\vec{r} = r \vec{\hat{e}}_r$}
\begin{subequations}
\begin{flalign}
&\mathcal{R}^{-}    :\hspace{4ex} \vec{\hat{e}}_{-}^{*}     \cdot \left( \vec{r} \, \Psi_{m,l,k}^{0,\alpha} \right) = \Psi_{m,l,\bigcdot}^{-,\alpha} \cdot \frac{1}{2} \mathcal{I}_{b}^{\dagger} & \\
&\mathcal{R}^{+}    :\hspace{4ex} \vec{\hat{e}}_{+}^{*}     \cdot \left( \vec{r} \, \Psi_{m,l,k}^{0,\alpha} \right) = \Psi_{m,l,\bigcdot}^{+,\alpha} \cdot \frac{1}{2} \mathcal{I}_{b} & \\
&\mathcal{R}^{\zero}:\hspace{4ex} \vec{\hat{e}}_{\zero}^{*} \cdot \left( \vec{r} \, \Psi_{m,l,k}^{0,\alpha} \right) = \Psi_{m,l+1,\bigcdot}^{0,\alpha} \cdot \frac{1}{\sqrt{2}} c_{l}^{(\alpha)} \mathcal{I}_{a} + \Psi_{m,l-1,\bigcdot}^{0,\alpha} \cdot \frac{1}{\sqrt{2}} d_{l}^{(\alpha)} \mathcal{I}_{a}^{\dagger} &
\end{flalign}
\end{subequations}

\subsection{Spherical Radial Component}
\begin{subequations}
\begin{flalign}
&\mathcal{E}^{+}    :\hspace{4ex} \vec{r} \cdot \left( \vec{\hat{e}}_{-}     \Psi_{m,l,k}^{-,\alpha} \right) = \Psi_{m,l,\bigcdot}^{0,\alpha} \cdot \frac{1}{2} \mathcal{I}_{b} & \\
&\mathcal{E}^{-}    :\hspace{4ex} \vec{r} \cdot \left( \vec{\hat{e}}_{+}     \Psi_{m,l,k}^{+,\alpha} \right) = \Psi_{m,l,\bigcdot}^{0,\alpha} \cdot \frac{1}{2} \mathcal{I}_{b}^{\dagger} & \\
&\mathcal{E}^{\zero}:\hspace{4ex} \vec{r} \cdot \left( \vec{\hat{e}}_{\zero} \Psi_{m,l,k}^{0,\alpha} \right) = \Psi_{m,l+1,\bigcdot}^{0,\alpha} \cdot \frac{1}{\sqrt{2}} c_{l}^{(\alpha)} \mathcal{I}_{a} + \Psi_{m,l-1,\bigcdot}^{0,\alpha} \cdot \frac{1}{\sqrt{2}} d_{l}^{(\alpha)} \mathcal{I}_{a}^{\dagger} &
\end{flalign}
\end{subequations}

\subsection{Multiplication by $1-r^2$}
First,
\begin{equation}
1 - r^2 = (1-\eta^2)(1-s^2) = \frac{1}{2} (1-\eta^2) (1 - t).
\end{equation}
Since
\begin{equation}
(1 - \eta^2) P_{l}^{(\alpha,\alpha)}(\eta) = \left( \mathcal{I}_{a} \mathcal{I}_{b} \right)^{\dagger} P_{l}^{(\alpha,\alpha)}(\eta) = \gamma_{l}^{(\alpha-1)} P_{l}^{(\alpha-1,\alpha-1)}(\eta) - \delta_{l}^{(\alpha-1)} P_{l+2}^{(\alpha-1,\alpha-1)}(\eta),
\end{equation}
we have
\begin{flalign}
&\mathcal{S}:\hspace{4ex} \left( 1 - r^2 \right) \Psi_{m,l,k}^{\sigma,\alpha} = \Psi_{m,l,\bigcdot}^{\sigma,\alpha-1} \cdot \left( +\frac{1}{2} \gamma_{l}^{(\alpha-1)} \mathcal{I}_{a}^{\dagger} \right) + \Psi_{m,l+2,\bigcdot}^{\sigma,\alpha-1} \cdot \left( -\frac{1}{2} \delta_{l}^{(\alpha-1)} \mathcal{I}_{a} \right) &
\end{flalign}

\subsection{Conversion}
\begin{flalign}
&\mathcal{I}_{\alpha}:\hspace{4ex} \Psi_{m,l,k}^{\sigma,\alpha} = \Psi_{m,l,\bigcdot}^{\sigma,\alpha+1} \cdot \left( +\gamma_{l}^{(\alpha)} \mathcal{I}_{a} \right) + \Psi_{m,l-2,\bigcdot}^{\sigma,\alpha+1} \cdot \left( - \delta_{l}^{(\alpha)} \mathcal{I}_{a}^{\dagger} \right) &
\end{flalign}

\subsection{Boundary Evaluation}
Our expansion for a single azimuthal mode evaluated on the surface of the ball, $\eta = \pm 1$, takes the form
\begin{equation}
f(t, \phi, \eta = \pm 1) = e^{i m \phi} (1+t)^{\frac{m}{2}} \sum_{l=0}^L (1-t)^{\frac{l}{2}} P_{l}^{(\alpha,\alpha)}(\pm 1) \sum_{k} \widehat{F}_{m,l,k} P_{k}^{(l+\alpha+\frac{1}{2},m)}(t).
\end{equation}
Since the basis decouples in the $\phi$ direction we must satisfy boundary data for each $m$ independently.  We thus derive the boundary evaluation operator
for fixed $m$ and drop the first index for notational convenience, so that ${ \widehat{F}_{m,l,k} \mapsto \widehat{F}_{l,k} }$.
In order to find a relationship between the expansion coefficients $\widehat{F}_{m,l,k}$ and the value of the expansion on the boundary we need to remove
the $t$ dependence.  We note that the $a$-lowering operator $\mathcal{I}_{a}^{\dagger}$ is equivalent to multiplication by $(1-t)$:
\begin{equation}
(1-t) P_{k}^{(a,b)}(t) = \mathcal{I}_{a}^{\dagger} P_{k}^{(a,b)}(t) \in \mathcal{H}(a-1,b),
\end{equation}
where $\mathcal{H}(a,b)$ here denotes the Hilbert space induced by the Jacobi integral weight $(1-t)^{a} (1+t)^{b}$.
By splitting into the even and odd $l$ indices and repeatedly applying $\mathcal{I}_{a}^{\dagger}$ then $\mathcal{I}_{a}$ we find:
\begin{equation}
\begin{aligned}
f(t, \phi, \pm 1) \propto 
    & (1-t)^{\frac{0}{2}} P_0(\pm 1) \sum_k \widehat{F}_{0,k} P_{k}^{(0+\frac{1}{2},m)}(t) + \hdots + 
      (1-t)^{\frac{L}{2}} P_L(\pm 1) \sum_k \widehat{F}_{L,k} P_{k}^{(L+\frac{1}{2},m)}(t) + \\
    & (1-t)^{\frac{1}{2}} P_1(\pm 1) \sum_k \widehat{F}_{1,k} P_{k}^{(1+\frac{1}{2},m)}(t) + \hdots + 
      (1-t)^{\frac{L-1}{2}} P_{L-1}(\pm 1) \sum_k \widehat{F}_{L-1,k} P_{k}^{(L-1+\frac{1}{2},m)}(t) \\
= 
    & \hphantom{(1-t)^{\frac{1}{2}}} \bigg( P_0(\pm 1) \sum_k \widehat{F}_{0,k} \left[ \left(\mathcal{I}_{a}\right)^{\frac{L}{2}} P_{\bigcdot}^{(0+\frac{1}{2},m)}(t) \right]_{k} \\
    &\hspace{10ex} + \hdots + 
      P_L(\pm 1) \sum_k \widehat{F}_{L,k} \left[ \left(\mathcal{I}_{a}^{\dagger}\right)^{\frac{L}{2}} P_{\bigcdot}^{(L+\frac{1}{2},m)}(t) \right]_{k} \bigg) + \\
    & (1-t)^{\frac{1}{2}} \bigg( 
      P_1(\pm 1) \sum_k \widehat{F}_{1,k} \left[ \left(\mathcal{I}_{a}\right)^{\frac{L-2}{2}} P_{\bigcdot}^{(1+\frac{1}{2},m)}(t) \right]_{k} \\
    &\hspace{10ex} + \hdots + P_{L-1}(\pm 1) \sum_k \widehat{F}_{L-1,k} \left[ \left(\mathcal{I}_{a}^{\dagger}\right)^{\frac{L-2}{2}} P_{\bigcdot}^{(L-1+\frac{1}{2},m)}(t) \right]_{k} \bigg).
\end{aligned}
\end{equation}
To parse the above result note the even $l$ terms take the form
\begin{equation}
\begin{aligned}
& P_0(\pm 1) \sum_k \widehat{F}_{0,k} \left[ \left(\mathcal{I}_{a}\right)^{\frac{L}{2}} P_{\bigcdot}^{(0+\frac{1}{2},m)}(t) \right]_{k}  + 
  P_2(\pm 1) \sum_k \widehat{F}_{2,k} \left[ \left(\mathcal{I}_{a}\right)^{\frac{L-2}{2}} \mathcal{I}_{a}^{\dagger} P_{\bigcdot}^{(2+\frac{1}{2},m)}(t) \right]_{k} + \hdots + \\
& P_{L-2}(\pm 1) \sum_k \widehat{F}_{L-2,k} \left[ \mathcal{I}_{a} \left(\mathcal{I}_{a}^{\dagger}\right)^{\frac{L-2}{2}} P_{\bigcdot}^{(L-2+\frac{1}{2},m)}(t) \right]_{k} +
  P_L(\pm 1) \sum_k \widehat{F}_{L,k} \left[ \left(\mathcal{I}_{a}^{\dagger}\right)^{\frac{L}{2}} P_{\bigcdot}^{(L+\frac{1}{2},m)}(t) \right]_{k} .
\end{aligned}
\end{equation}
Each successive term gets one more lowering $\mathcal{I}_{a}^\dagger$ operator to account for the additional factor of $(1-t$) and one fewer $\mathcal{I}_{a}$
operator since the $a$ index starts closer to the final $a$ index $\frac{L}{2}$ common to all even terms.  We end up in the Jacobi polynomial
space $\mathcal{H}\left(\frac{L}{2}, m\right)$ for even $l$ and $\mathcal{H}\left(\frac{L-2}{2},m\right)$ for odd $l$.  This enables us to reverse the
summations over $l$ and $k$ and we end up with
\begin{equation}
\begin{aligned}
f(t, \phi, \pm 1) \propto 
    \hspace{3ex} &
    \sum_k P_{k}^{\left(\frac{L}{2},m\right)}(t) \left( P_0(\pm 1) \left[ \left(\mathcal{I}_{a}\right)^{\frac{L}{2}} \widehat{F}_{0,\bigcdot} \right]_{k} + \hdots 
        + P_L(\pm 1) \left[ \left(\mathcal{I}_{a}^{\dagger}\right)^{\frac{L}{2}} \widehat{F}_{L,\bigcdot} \right]_{k} \right) \\
   &\hspace{-13ex} + (1-t)^{\frac{1}{2}}  \left[ 
        \sum_k P_{k}^{\left(\frac{L-2}{2},m\right)}(t) \left( P_1(\pm 1) \left[ \left(\mathcal{I}_{a}\right)^{\frac{L-2}{2}} \widehat{F}_{1,\bigcdot} \right]_{k} + \hdots 
        + P_{L-1}(\pm 1) \left[ \left(\mathcal{I}_{a}^{\dagger}\right)^{\frac{L-2}{2}} \widehat{F}_{L-1,\bigcdot} \right]_{k} \right)
    \right].
\end{aligned}
\end{equation}
Hence to force a field to zero on the boundary we require, for each $k$,
\begin{equation}
\begin{aligned}
0 &= P_{0}(\pm 1) \left[ \left(\mathcal{I}_{a}\right)^{\frac{L}{2}} \widehat{F}_{0,\bigcdot} \right]_{k} 
    + P_{2}(\pm 1) \left[ \left(\mathcal{I}_{a}\right)^{\frac{L}{2}-1} \mathcal{I}_{a}^{\dagger} \widehat{F}_{2,\bigcdot} \right]_{k} + \hdots \\
   &\hspace{4ex}\hdots + P_{L-2}(\pm 1) \left[ \mathcal{I}_{a} \left(\mathcal{I}_{a}^{\dagger}\right)^{\frac{L}{2}-1} \widehat{F}_{L-2,\bigcdot} \right]_{k}
    + P_{L}(\pm 1) \left[ \left(\mathcal{I}_{a}^{\dagger}\right)^{\frac{L}{2}} \widehat{F}_{L,\bigcdot} \right]_{k}
\end{aligned}
\end{equation}
and
\begin{equation}
\begin{aligned}
0 &= P_{1}(\pm 1) \left[ \left(\mathcal{I}_{a}\right)^{\frac{L}{2}-1} \widehat{F}_{1,\bigcdot} \right]_{k} 
    + P_{3}(\pm 1) \left[ \left(\mathcal{I}_{a}\right)^{\frac{L}{2}-2} \mathcal{I}_{a}^{\dagger} \widehat{F}_{3,\bigcdot} \right]_{k} + \hdots  \\
   &\hspace{4ex}\hdots + P_{L-3}(\pm 1) \left[ \mathcal{I}_{a} \left(\mathcal{I}_{a}^{\dagger}\right)^{\frac{L}{2}-2} \widehat{F}_{3,\bigcdot} \right]_{k} + \hdots 
    + P_{L-1}(\pm 1) \left[ \left(\mathcal{I}_{a}^{\dagger}\right)^{\frac{L}{2}-1} \widehat{F}_{L-1,\bigcdot} \right]_{k}.
\end{aligned}
\end{equation}
Note that this derivation assumes $L$ is even.  Counting the number of raising and lowering operations changes slightly in the case $L$ is odd.

\begin{figure}[H]
\centering
\includegraphics[width=0.8\linewidth]{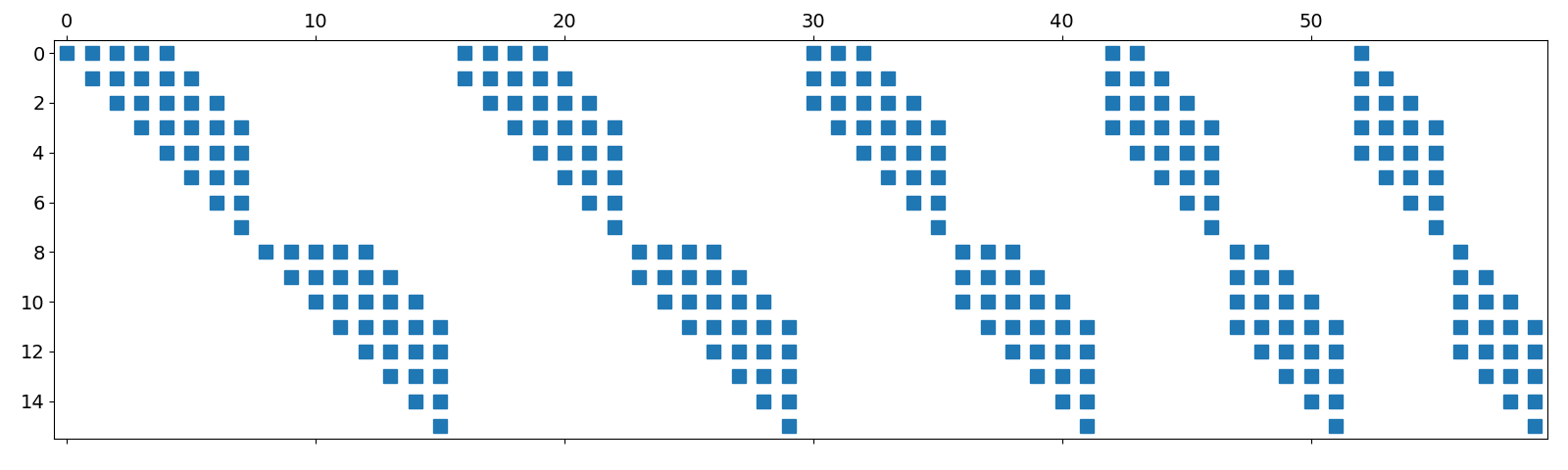}
\caption{Sparsity structure of the boundary evaluation operator for $\left( L_{\text{max}}, N_{\text{max}} \right) = \left( 10, 8 \right)$.}
\label{fig:spy-boundary}
\end{figure}

Figure~\ref{fig:spy-boundary} plots the sparsity structure of the boundary evaluation operator.  The first set of rows are the coupling of the even $l$ modes
while the second set of rows are the coupling of the odd $l$ modes.  The sub-blocks reduce in number of columns corresponding to the triangular truncation employed with the basis,
namely $N(l) = N_{\text{max}} - \floor{\frac{l}{2}}$.

\section{Regularity at the Equator}
We apply Frobenius's method for the indicial exponent to Laplace's equation ${ \nabla^2 f = 0 }$ in the ball to obtain leading order behavior at the equator.
From the result in the disk we know scalars with azimuthal wavenumber $m$ must behave like
\begin{equation}
f(s) \sim e^{i m \phi} s^{\left| m \right|} F(s^2),
\end{equation}
where $F(s^2)$ is even and analytic in a neighborhood of ${ s = 0 }$.  We seek the leading order behavior of $F$ at
${ s = 1 }$.  We therefore expand in a power series there assuming a vertical polynomial of degree $l$, maintaining
the requirement that $F$ be even in $s$:
\begin{equation}
f = e^{i m \phi} s^{\left| m \right|} \eta^{l} \sum_{k = 0}^{\infty} f_{k} (1 - s^2)^{k+a},
\end{equation}
where $a$ is the yet determined leading order behavior at the equator.  We first note the Laplacian in spherindrical coordinates takes the form
\begin{equation}
\begin{aligned}
\nabla^2 f 
   = \frac{1}{s} \partial_{s} \left( s \partial_{s} f \right) + \frac{1}{s^2} \partial_{\phi \phi} f + \frac{2 s \eta}{1 - s^2} \partial_{s} \partial_{\eta} f + \frac{1 - s^2 + s^2 \eta^2}{(1-s^2)^2} \partial_{\eta \eta} f + \frac{2 +
   s^2}{(1-s^2)^2} \eta \partial_{\eta} f.
\end{aligned}
\end{equation}
Plugging in our expansion we find leading order behavior
\begin{equation}
\begin{aligned}
\nabla^2 f &\sim e^{i m \phi} s^{\left| m \right|} (1-s^2)^{a-2} \eta^{l - 2} \\
    &\hspace{8ex} \times \left\{ l (l-1) \left(1-s^2\right) + (l-2a)\left[2+2\left(1-s^2\right) \left| m \right| +(l-2a)s^2\right]\eta^2 \right\} f_{0}.
\end{aligned}
\end{equation}
In a neighborhood of ${ s = 1 }$ we have to leading order
\begin{equation}
\begin{aligned}
\nabla^2 f &\sim e^{i m \phi} (1-s^2)^{a-2} \eta^{l} \times \left\{ (l-2a)\left[l+2-2a\right] \right\} f_{0}, \hspace{4ex} s \to 1.
\end{aligned}
\end{equation}
This expansion has two nontrivial solutions (i.e. ${ f_{0} \ne 0 }$) to $\nabla^2 f = 0$:
\begin{equation}
a = \frac{l}{2}, \hspace{2ex} \frac{l+2}{2}.
\end{equation}
We discard the solution ${ a = \frac{l}{2} + 1 }$ since it cannot represent fields constant throughout the ball.
Thus for a field to satisfy Laplace's equation in the stretched sphere it must have leading behavior
\begin{equation}
f \sim e^{i m \phi} s^{\left| m \right|} \eta^{l} (1-s^2)^{\frac{l}{2}},
\end{equation}
which is simply the Cartesian polynomial ${ \left( x + i y \right)^{m} z^{l} }$.
This is the motivation for the ${ (1-t)^{\frac{l}{2}} }$ prefactor in our basis functions $\Psi_{m,l,k}^{\sigma, \alpha}$ defined in (\ref{eqn:spherinder_basis}).

\end{appendices}

\bibliographystyle{unsrt}
\bibliography{References}

\end{document}